\documentclass[12pt]{scrartcl}

\usepackage{setspace}

\usepackage{color}
\usepackage{array}
\usepackage{supertabular}
\usepackage{hhline}
\usepackage{multirow}
\usepackage{multicol}
\usepackage{textgreek}
\usepackage{upgreek}

\usepackage{wrapfig}

 \usepackage{mathtools}

\usepackage{stmaryrd}

\usepackage[listings,theorems]{tcolorbox}

\usepackage{pifont}

\usepackage[thicklines]{cancel}

\usepackage[LGR,T1]{fontenc}

\usepackage{tikz}

\usetikzlibrary{hobby}

\usepackage[customcolors]{hf-tikz}

\usepackage{circuitikz}

\usetikzlibrary{decorations.pathreplacing,arrows.meta}

\usetikzlibrary{calc,patterns,decorations.markings}
\usetikzlibrary{fadings}
\usetikzlibrary{angles,quotes}
\usepackage{tikz-3dplot}

\hfsetfillcolor{yellow!10}
\hfsetbordercolor{darkblue}

\def\uphrulefill{%
  \textcolor{darkblue!70}{\leavevmode\leaders\hrule height 1.5pt\hfill\kern0pt \rule[-0.56ex]{1.5pt}{5pt}}\\}
  
  \def\downhrulefill{%
  \textcolor{darkblue!70}{\leavevmode\leaders\hrule height 1.5pt\hfill\kern0pt \rule[0.0ex]{1.5pt}{5pt}}\\}

\definecolor{darkgreen}{rgb}{0.0, 0.2, 0.13}
\definecolor{darkred}{rgb}{0.55, 0.0, 0.0}
\definecolor{darkgreen}{rgb}{0.0, 0.2, 0.13}
\definecolor{darkmagenta}{rgb}{0.55, 0.0, 0.55}

\usepackage{pgfplots}

\pgfplotsset{compat=1.10}
\usepgfplotslibrary{fillbetween}
\usetikzlibrary{patterns}

 \usetikzlibrary{patterns}

\usepackage[latin1]{inputenc}

\definecolor{darkspringgreen}{rgb}{0.09, 0.45, 0.27}
\definecolor{oxfordblue}{rgb}{0.0, 0.13, 0.28}
\definecolor{upmaroon}{rgb}{0.48, 0.07, 0.07}
\definecolor{darkmagenta}{rgb}{0.55, 0.0, 0.55}
\definecolor{byzantium}{rgb}{0.44, 0.16, 0.39}
\definecolor{firebrick}{rgb}{0.7, 0.13, 0.13}
\definecolor{jazzberryjam}{rgb}{0.65, 0.04, 0.37}
\definecolor{auburn}{rgb}{0.43, 0.21, 0.1}
\definecolor{champagne}{rgb}{0.97, 0.91, 0.81}
\definecolor{mediumaquamarine}{rgb}{0.4, 0.8, 0.67}
\definecolor{teadarkspringgreen}{rgb}{0.82, 0.94, 0.75}
\definecolor{uclagold}{rgb}{1.0, 0.7, 0.0} 
\definecolor{indigo}{rgb}{0.29, 0.0, 0.51}
\definecolor{anti-flashwhite}{rgb}{0.95, 0.95, 0.96}
\definecolor{darkblue}{rgb}{0.0, 0.0, 0.55}

\usepackage{graphicx}
\usepackage{fancybox}

\usepackage[latin1]{inputenc}

\usepackage{tcolorbox}

\def\uphrulefill{%
  \textcolor{darkblue!70}{\leavevmode\leaders\hrule height 1.5pt\hfill\kern0pt \rule[-0.56ex]{1.5pt}{5pt}}\\}
  
  \def\downhrulefill{%
  \textcolor{darkblue!70}{\leavevmode\leaders\hrule height 1.5pt\hfill\kern0pt \rule[0.0ex]{1.5pt}{5pt}}\\}

 \usepackage{multirow}

\usepackage{graphicx,fancybox,manfnt,wrapfig,multicol,mathrsfs}

\usepackage{amsmath,amssymb,color,makeidx,rotating,pifont,amsthm}

\usepackage{mathtools}

\DeclareMathAlphabet{\boondoxcal}{U}{BOONDOX-cal}{m}{n}

\newcommand{\mys}{\boondoxcal{s}}

\mathtoolsset{showonlyrefs=true}  

\newtheorem{prop}{Proposition}
\newtheorem{defin}{Definition}

\allowdisplaybreaks

\begin{document}

\begin{center}
\LARGE
\bfseries\sffamily A Generalization of the Fox $\textbf{H}$-function 
\end{center}
\normalfont
\normalsize

\begin{center}
{\bfseries\sffamily \large Jayme Vaz}$^{a,b,}$\footnote{\texttt{vaz@unicamp.br}}
\end{center}

\normalsize\rmfamily\sffamily

\begin{center}
\begin{minipage}[c]{14cm}
\begin{flushleft}
$^{a}$ DSS - Sapienza Universit\`a di Roma, Piazzale Aldo Moro 5, Rome 00185, Italy \\
$^{b}$ IMECC - Universidade Estadual de Campinas, 13087-859 Campinas, SP, Brazil
\end{flushleft}
\end{minipage}
\end{center}

\bigskip

\small\rmfamily

\begin{center}
\begin{minipage}[c]{12cm}
\centerline{\bfseries\sffamily Abstract}
In this paper we present a generalization of the Fox $H$-function called Fox-Barnes $J$-function. Like the Fox $H$-function, it is defined as a contour integral in the complex plane, but instead of an integrand given by a ratio of products of gamma functions involving several parameters, we use a ratio of products of double gamma functions. We study the conditions for its existence and how to choose a contour of integration based on the involved parameters. We discuss how the Fox $H$-function appears as a particular case and prove some properties of the Fox-Barnes $J$-function. As an application, we show how the Laplace transform of the Kilbas-Saigo function can be conveniently written in terms of the Fox-Barnes $J$-function, even in cases where the usual series representation is not convergent. 
\end{minipage}
\end{center}

\medskip

\section{Introduction}

The generalization of special functions plays a important role in mathematics, 
enabling the unification of different results, expanding the scope of 
analytical methods and the emergence of new applications. A well-known  
example of this is the hypergeometric function, whose structure 
encompasses a wide range of classical functions, 
making it a powerful tool for solving differential equations and 
treating special integrals \cite{HG1,HG2,HG3}. 

Advancing this process of generalization, Meijer $G$-function emerges 
as an even broader extension of the hypergeometric function \cite{Meijer1,Meijer2}. 
Its definition, based on contour integrals in the complex plane (Mellin-Barnes integral), 
gives it great analytical and computational flexibility.
The Fox $H$-function represents a new step in this process of abstraction and generalization 
\cite{Fox1,Fox2} as 
it generalizes the Meijer $G$-function by incorporating a richer structure 
of parameters at the poles of the gamma functions that comprise its definition. 
Thus, Fox $H$-function became one of the most general objects within the theory of 
special functions defined by Mellin-Barnes integrals.

These functions (hypergeometric, Meijer $G$ and Fox $H$) have been widely used 
in statistics, particularly in modeling probability distributions and solving 
problems involving integral transforms. For example, the hypergeometric function 
appears naturally in the hypergeometric distribution and in models involving 
generalized moments \cite{HGApp}. Meijer $G$-function is used to express densities 
of distributions such as the lognormal, the generalized variance distribution, 
and various alpha-stable distributions. Fox $H$-function  
finds application in describing distributions with heavy tails, such as those 
that arise in non-Gaussian stochastic processes, 
in memory-based waiting time models, and in anomalous diffusion phenomena \cite{BCS,BCS2,BoCS,Denisov,Ranga}.

Given the breadth and depth of applications of these functions in statistics 
and related fields, it is natural to seek even more general forms of 
functional representation, in the hope of modelling more complex behaviours and 
introducing new types of parametric structures.  In fact, recent works \cite{Nikeghbali, Ostrovsky, Simon0, Simon2, Kusnetsov} have shown that a broad 
class of probability laws, ranging from Barnes beta distributions and 
infinite products of beta random variables to sums and products of 
generalized gamma convolution variables, are naturally described through
their Mellin transforms and Barnes' double gamma function \cite{Genesis}. In this 
context, a generalization of the Fox $H$-function to encompass the double gamma
function may provide an integral representation for probability densities
whose Mellin transforms involve ratios of double gamma functions, extending therefore the classical Fox $H$-function framework to settings motivated by modern probabilistic theory and stochastic process modelling. 

Thus in this work we discuss a possible generalization of Fox $H$-function, which we will call the Fox-Barnes $J$-function, in analogy 
with the generalization of the Meijer $G$-function proposed by Karp and 
Kuznetsov in \cite{Kusnetsov}. The idea is to still use a Mellin-Barnes integral, but instead of the gamma function, to use the double gamma function. Certainly calling such a function as $I$-function would be a good name for a generalization of a $H$-function. However, there are other proposals for functions that generalize Fox $H$-function, as in \cite{I-2,I-1}, where powers of gamma functions are used in the integrand instead of gamma functions, or as in \cite{I-3} where a sum of products of gamma functions is used in the denominator of the 
integrand, and all these functions were 
already called $I$-functions. In order to avoid confusion with these other proposals, we chose the next letter of the English alphabet for our proposal. 
Moreover the explicit mention of Barnes' name suggests the use of the double gamma function in place of the gamma function.

We organized this paper as follows: in Section~\ref{sec.2} we review the definition and main properties of the double gamma function, in Section~\ref{sec.3} we define the Fox-Barnes $J$-function and discuss conditions for its existence, in Section~\ref{sec.4} we show how the Fox $H$-function appears as a particular case of the Fox-Barnes $J$-function, in Section~\ref{sec.5} we discuss some of its properties and in Section~\ref{sec.6} we discuss an application of it in the calculation of the Laplace transform of the Kilbas-Saigo function for a case where a traditional series representation is not convergent.  The motivation for this problem comes from the studies about 
stochastic models governed by stretched fractional dynamics \cite{BLV},  
where in a generalized renewal processes modelled by the Kilbas-Saigo function, the parameter regime that ensures finite expected value interarrival times leads to a formally divergent Laplace series. This makes it necessary to go beyond classical series-based transforms and motivates the use of Fox-Barnes $J$ functions, which provide a well-defined representation of the relevant statistical structure.

\section{The double gamma function}
  \label{sec.2}

The double gamma function, 
denoted by $G(z; \tau)$, with $z\in \mathbb{C}$ and $\tau > 0$,
is defined as \cite{Alexanian,SEFC2nd,Vaz_Dim}
\begin{equation}
\label{ap.B.G.tau.1}
\begin{split}
G(z;\tau) = & \frac{z}{\tau} \exp\left[a(\tau)\frac{z}{\tau} + b(\tau)\frac{z^2}{2\tau}\right] \\
& \cdot \prod_{m,n\in \mathbb{N}_\ast^2} \left[
\left(1+\frac{z}{m\tau+n}\right) \exp\left(-\frac{z}{m\tau+n}+\frac{z^2}{2(m\tau+n)^2}\right)\right] ,
\end{split}
\end{equation}
where $\mathbb{N}_\ast^2 = \mathbb{N}\setminus \{(0,0)\}$, and $a(\tau)$, $b(\tau)$ are given by
\begin{equation}
\begin{split}
& a(\tau) = \gamma \tau + \frac{\tau}{2}\log{(2\pi \tau)} + \frac{1}{2}\log\tau - \tau C(\tau) , \\
& b(\tau) = - \frac{\pi^2 \tau^2}{6} - \tau \log\tau - \tau^2 D(\tau) ,
\end{split}
\end{equation}
where $\gamma$ is the Euler-Mascheroni constant, and
\begin{align}
& C(\tau) = \lim_{m\to \infty} \left[\sum_{k=1}^{m-1}\psi(k\tau) + \frac{1}{2}\psi(m\tau) - \frac{1}{\tau}
\log\left(\frac{\Gamma(m\tau)}{\sqrt{2\pi}}\right) \right] , \\
& D(\tau) = \lim_{m\to \infty} \left[\sum_{k=1}^{m-1}\psi^\prime(k\tau) + \frac{1}{2}\psi^\prime(m\tau) - \frac{1}{\tau}\psi(m\tau) \right] ,
\end{align}
with $\psi(\cdot)$ being the digamma function.

An alternative definition of double gamma function is
\begin{equation}
\label{ap.B.G.tau.2}
G(z;\tau) = \frac{1}{\tau\Gamma(z)} \exp\left[\tilde{a}(\tau)\frac{z}{\tau} + \tilde{b}(\tau)
\frac{z^2}{2\tau^2}\right] \prod_{m=1}^\infty \frac{\Gamma(m\tau)}{\Gamma(z+m\tau)}
\exp\left[z\psi(m\tau)+ \frac{z^2}{2}\psi^\prime(m\tau) \right] ,
\end{equation}
where
\begin{equation}
\tilde{a}(\tau) = a(\tau) - \gamma\tau , \qquad \tilde{b}(\tau) = b(\tau) + \frac{\pi^2 \tau^2}{6} .
\end{equation}
The Barnes $G$-function $G(z)$ is a particular case of the double
gamma function $G(z;\tau)$ for $\tau = 1$. Note that $C(1) = 1/2$ and $D(1) = 1+ \gamma$, as shown in \cite{Genesis}.

The double gamma function is an entire function and its zeros are located
at $z = z_{mn} = -m\tau-n$ ($m,n = 0,1,2,\ldots$).  
It is defined in such a way that \cite{Genesis}
\begin{equation}
\label{ap.B.G.tau.(1)}
G(1;\tau) = 1 .
\end{equation}
It satisfies the functional relations
\begin{equation}
\label{ap.B.general.G.tau}
G(z+1;\tau) = \Gamma(z/\tau) G(z;\tau)
\end{equation}
and
\begin{equation}
\label{ap.B.general.G.tau.2}
G(z+\tau;\tau) = (2\pi)^{\frac{\tau-1}{2}} \tau^{-z+\frac{1}{2}} \Gamma(z) G(z;\tau) .
\end{equation}
Using Eq.\eqref{ap.B.general.G.tau} recursively, we obtain
\begin{equation}
\label{ap.B.general.G.tau.n}
G(z+k;\tau) = G(z,\tau) \prod_{j=0}^{k-1} \Gamma[(z+j)/\tau] .
\end{equation}
Another interesting property of $G(z;\tau)$ is the modular transformation 
\begin{equation}
\label{modular}
G(z;\tau) = (2\pi)^{\frac{z}{2\tau}(\tau-1)}\,\tau^{\frac{z-z^2}{2\tau}+\frac{z}{2}-1} G\left(\frac{z}{\tau};
\frac{1}{\tau}\right) . 
\end{equation}
Indeed it is not difficult to see that eq.\eqref{ap.B.general.G.tau} and eq.\eqref{ap.B.general.G.tau.2} 
are related by eq.\eqref{modular}.

For large $z$ with $|\operatorname{arg} z| < \pi$, it holds \cite{Alexanian,Spreafico} 
\begin{equation}
\label{stirling.G}
\log G(z;\tau) = [a_2(\tau) z^2 + a_1(\tau) z + a_0(\tau)]\log{z} + b_2(\tau)z^2 + b_1(\tau)z + \mathcal{O}(1) 
\end{equation}
where 
\begin{align}
& a_2(\tau) = \frac{1}{2\tau} , \qquad a_1(\tau) = -\frac{1}{2}\left(1+\frac{1}{\tau}\right) , \qquad 
a_0(\tau) =\frac{\tau}{12} + \frac{1}{4} + \frac{1}{12\tau} , \\[1ex]
& b_2(\tau) = -\frac{1}{2\tau}\left(\frac{3}{2}+\log\tau\right) , \qquad 
b_1(\tau) = \frac{1}{2}\left[\left(1+\frac{1}{\tau}\right)(1+\log\tau) + \log{2\pi}\right] . 
\end{align}

An integral representation for $G(z;\tau)$ was provided in \cite{lawrie}, that is,
\begin{equation}
\label{double.gamma.int.rep}
\begin{aligned}
\log G(z;\tau) & = \int_0^1 \bigg[ \frac{r^{z-1}}{(r-1)(r^\tau-1)} -
\frac{z^2}{2\tau} r^{\tau-1} - z r^{\tau-1}\left(\frac{2-r^\tau}{r^\tau-1}-\frac{1}{2\tau}\right) \\
& -r^{\tau-1} + \frac{1}{r-1} - \frac{r^{\tau-1}}{(r-1)(r^\tau-1)}\bigg] \frac{dr}{\ln{r}} .
\end{aligned}
\end{equation}
It is convergent for $\operatorname{Re}z > 0$ and $\delta > 0$.

\section{The Fox-Barnes $\boldsymbol{J}$-function}
 \label{sec.3}

We define the Fox-Barnes $J$-function as a generalization of the Meijer-Barnes $K$-function \cite{Kusnetsov} in 
the same way the Fox $H$-function is defined as a generalization of the Meijer $G$-function, that is, 
\begin{equation}
\label{fox-barnes}
J^{m,n}_{p,q}\left[ z;\tau,\varepsilon \, \bigg| \, \begin{matrix} (a_1,\alpha_1),\ldots,(a_p,\alpha_p) \\ (b_1,\beta_1),\ldots,(b_q,\beta_q)\end{matrix}\right] = 
\frac{1}{2\pi i} \int_{\mathcal{C}} \mathcal{J}^{m,n}_{p,q}(s) \operatorname{e}^{\pi \varepsilon s^2/\tau} z^{-s}\, ds ,
\end{equation}
with $z \in \mathbb{C}$, $\tau > 0$ and $\mathcal{J}^{m,n}_{p,q}(s)$ given by 
\begin{equation}
\label{fox-barnes.integrand}
\mathcal{J}^{m,n}_{p,q}(s) = \frac{\prod_{i=1}^m G(b_i+\beta_i s;\tau) \cdot \prod_{i=1}^n G(1+\tau -a_i-\alpha_i s;\tau)}{ 
\prod_{i=m+1}^q G(1+\tau -b_i-\beta_i s;\tau)\cdot \prod_{i=n+1}^p G(a_i +\alpha_i s;\tau) }   
\end{equation}
where $\alpha_i,\beta_i > 0$, $\tau > 0$, $\varepsilon \in \mathbb{C}$ and $a_i, b_i \in \mathbb{C}$. The parameter $\varepsilon$ will be called \textit{compensation parameter} and the reason for its introduction will be clear later (see Proposition~\ref{prop.3} and Remark~2).  
When $\varepsilon$ is not present in the notation on the left side of eq.\eqref{fox-barnes} 
this means that $\varepsilon = 0$. 
A more compact notation we will use is 
\begin{equation}
J^{m,n}_{p,q}\left[ z;\tau,\varepsilon \, \bigg| \, \begin{matrix} (a_i,\alpha_i)_{1,p} \\ (b_i,\beta_i)_{1,q} \end{matrix}\right]  = 
J^{m,n}_{p,q}\left[ z;\tau,\varepsilon \, \bigg| \, \begin{matrix} (a_1,\alpha_1),\ldots,(a_p,\alpha_p) \\ (b_1,\beta_1),\ldots,(b_q,\beta_q)\end{matrix}\right]  .
\end{equation}

\medskip

 The contour of integration $\mathcal{C}$ separates the poles of $\mathcal{J}^{m,n}_{p,q}(s)$ 
due to the zeros of $G(1+\tau -b_i-\beta_i s;\tau)$, that is, 
\begin{equation}
\label{eq.poles.left}
\acute{\mys}^i_{m^\prime, n^\prime} = \frac{1}{\beta_i}\left(-b_i + m^\prime\tau + n^\prime\right), \qquad m^\prime,n^\prime = 1,2,\ldots ,\; i = m+1,\ldots, q,  
\end{equation}
 from the poles  of $\mathcal{J}^{m,n}_{p,q}(s)$ due to the zeros 
of $G(a_j +\alpha_j s;\tau)$, that is, 
\begin{equation}
\label{eq.poles.right}
\grave{\mys}^j_{m, n} = -\frac{1}{\alpha_j}\left(a_j + m \tau + n\right), \qquad m ,n  = 0,1,2,\ldots , \;
j = n+1,\ldots, p .
\end{equation}
We assume that the poles do not coincide, 
$$
\alpha_j (b_i - m^\prime\tau-n^\prime) \neq \beta_i(a_j + m\tau +n)   
$$
for $i = m+1,\ldots, q$ and $j = n+1,\ldots, p$. 

Let $\mathcal{C}_{(\theta^-,\theta^+)}$ be a contour starting with a half-line 
in the lower half-plane with slope $\tan\theta^-$ and terminating 
with a half-line in the upper half-plane with slope $\tan\theta^+$, leaving to the right all 
poles $\acute{\mys}^i_{m^\prime,n^\prime}$ as 
in eq.\eqref{eq.poles.left} and to the left all the poles 
 $\grave{\mys}^j_{m,n}$ as in eq.\eqref{eq.poles.right} (see \cite{Kusnetsov} for some illustrations of these kind of contours). 
 We have a special interest in the contours  $\mathcal{C}_{ic\infty}$, 
 $\mathcal{C}_{-\infty}$ and $\mathcal{C}_{+\infty}$, 
since $\mathcal{C}_{ic\infty}$ is used in the inversion of the
Mellin and Laplace transforms, while $\mathcal{C}_{-\infty}$ and $\mathcal{C}_{+\infty}$ 
can be employed in conjunction with the Cauchy residue theorem. Let us recall that they are
defined as:  (i)  $\mathcal{C} = \mathcal{C}_{-\infty}$ is a left 
loop situated in a horizontal strip starting at $-\infty + i\varphi_1$ and
terminating at $-\infty + i \varphi_2$ with $-\infty < \varphi_1 < 0 <  \varphi_2 < \infty$ 
and encircling once in the positive direction all the poles $\grave{\mys}^j_{m,n}$ but none of the poles  $\acute{\mys}^i_{m^\prime,n^\prime}$ -- thus the contour $\mathcal{C}_{-\infty}$ corresponds to $\mathcal{C}_{(-\pi+\epsilon,\pi-\epsilon)}$ for $\epsilon \to 0$; (ii)  
 $\mathcal{C} = \mathcal{C}_{+\infty}$ is a right 
loop situated in a horizontal strip starting at $+\infty + i\varphi_2$ and
terminating at $+\infty + i \varphi_1$ with $-\infty < \varphi_1 < 0 < \varphi_2 < \infty$ 
and encircling once in the positive direction all the poles $\acute{\mys}^i_{m^\prime,n^\prime}$   but none of the poles  $\grave{\mys}^j_{m,n}$ -- thus 
the contour $\mathcal{C}_{+\infty}$ 
corresponds to $\mathcal{C}_{(0-\epsilon,0+\epsilon)}$ for $\epsilon \to 0$;  
(iii) $\mathcal{C} = \mathcal{C}_{ic\infty}$ is a contour starting at $c-i\infty$ 
and terminating at $c+i\infty$ ($c \in \mathbb{R}$) leaving to the right all 
poles $\acute{\mys}^i_{m^\prime,n^\prime}$  and to the left all the poles 
 $\grave{\mys}^j_{m,n}$ -- thus the contour  $\mathcal{C}_{ic\infty}$ corresponds 
 to $\mathcal{C}_{(-\frac{\pi}{2},\frac{\pi}{2})}$.
 
Let us define $\mathcal{A}$ and $\mathcal{B}$ as 
\begin{equation}
\label{def.A.B}
\mathcal{A} = -\underset{n+1\leq i\leq p}{\operatorname{min}} \;\frac{\operatorname{Re} a_i}{\alpha_i} , \qquad 
\mathcal{B} = \underset{m+1\leq i\leq q}{\operatorname{min}}\frac{1+\tau - \operatorname{Re} b_i}{\beta_i} .
\end{equation}
To simplify the problem, we will assume that $\mathcal{A} < \mathcal{B}$. Thus we have 
\begin{equation}
\operatorname{Re} \grave{\mys}^i_{m,n}  \leq \mathcal{A} < \mathcal{B} \leq \operatorname{Re} \acute{\mys}^i_{m^\prime,n^\prime}.
\end{equation}

In order to choose the contours, we need to analyse the behaviour of the
integrand in eq.\eqref{fox-barnes} for $|s| \to \infty$ along the 
above contours. 

\bigskip

\subsection{Existence}
\label{section.3.1}

Writing $s = c + R\operatorname{e}^{i\theta}$ with $c \in \mathbb{R}$ and recalling that 
$\sigma \in \mathbb{C}$ and $\nu, \tau > 0$, it follows from eq.\eqref{stirling.G}  by replacing $z$ by $\sigma \pm \nu s$, using $\operatorname{arg}(-s) = \theta - \pi$ for $\theta = 
\operatorname{arg}(s) > 0$ and $\operatorname{arg}(-s) = \theta + \pi$ for $\theta = \operatorname{arg}(s) < 0$, and collecting the terms involving $R$ and 
$\log{R}$,
that 
\begin{equation}
\begin{aligned}
 & \log|G(\sigma\pm \nu s;\tau)| = \left( A_2(\tau,\nu)\cos{2\theta}\right) R^2\log{R} 
+ [B_2(\tau,\nu)\cos{2\theta} - 
A_2(\tau,\nu) \, \theta_\pm \sin{2\theta} ] R^2 \\[1ex] 
& \qquad \pm \left[\operatorname{Re}(A_1(\tau,\sigma_\pm,\nu))\cos\theta - 
\operatorname{Im}(A_1(\tau,\sigma_\pm,\nu))\sin\theta\right] R\log{R}  \\[1ex]
& \qquad \pm \big[\operatorname{Re}(B_1(\tau,\sigma_\pm,\nu))\cos\theta-\operatorname{Re}(A_1(\tau,\sigma_\pm,\nu))\, \theta_\pm \sin\theta\\[1ex]
& \qquad - \operatorname{Im}(B_1(\tau,\sigma_\pm,\nu))\sin\theta - 
\operatorname{Im}(A_1(\tau,\sigma_\pm,\nu))\, \theta_\pm \cos\theta \big] R  
  +A_0(\tau,\sigma_\pm,\nu) \log{R}   
 + \mathcal{O}(1) , \\[1ex] 
\end{aligned}
\end{equation} 
where 
\begin{align*}
& \sigma_\pm = \sigma \pm c \nu , \qquad 
 A_2(\tau,\nu) = \nu^2 a_2(\tau) , \qquad 
 A_1(\tau,\sigma_\pm,\nu) =  \nu[2\sigma_\pm a_2(\tau) + a_1(\tau)] , \\[1ex]
& A_0(\tau,\sigma_\pm,\nu) =   \sigma_\pm^2 a_2(\tau) + \sigma_\pm a_1(\tau) + a_0(\tau) , \qquad 
 B_2(\tau, \nu) =  \nu^2[b_2(\tau) + a_2(\tau)\log\nu], \\[1ex]
& B_1(\tau,\sigma_\pm,\nu) =   \nu[2\sigma_\pm b_2(\tau) + b_1(\tau) + \sigma(1+2\log\nu) a_2(\tau) + \log{\nu}\, a_1(\tau)] ,\\[1ex]
& \theta_\pm  = \begin{cases} \theta, \; & \text{for $G(\sigma+\nu s;\tau)$} \\[1ex]
\theta + \big(1-2[\theta>0]\big)\pi = \begin{cases} \theta - \pi , \; \; & 
(\theta > 0) ,\\[1ex]
\theta + \pi , \; & (\theta \leq 0) , \end{cases}   \; \; & \text{for $G(\sigma-\nu s;\tau)$}
\end{cases}
\end{align*}
and $[P]$ is the Iverson bracket, i.e., 
$[P] = 1$ if $P$ is true and $[P] = 0$ otherwise. In \cite{BLV} and \cite{BLPV} similar expressions in other cases were studied.

Thus we can write for the integrand $\mathcal{I}(s)$ in eq.\eqref{fox-barnes} for $R \to \infty$ that 
\begin{equation}
\label{integrand}
\begin{split}
& \log|\mathcal{I}(s)| =  \log|\mathcal{J}^{m,n}_{p,q}(s)| + (-\cos\theta\log{|z|} + \sin\theta \operatorname{arg}z)R \\[1ex]
& + \frac{\pi}{\tau}\left( \operatorname{Re}\left(\varepsilon\right)\cos{2\theta}  -  \operatorname{Im}\left(\varepsilon\right)\sin{2\theta} \right) \, R^2  + 
\frac{2\pi c}{\tau}\left( \operatorname{Re}\left(\varepsilon\right)\cos{\theta}  -  \operatorname{Im}\left(\varepsilon\right)\sin{\theta} \right) \, R + \mathcal{O}(1)  
\end{split}
\end{equation}
with 
\begin{equation}
\begin{split}
\log|\mathcal{J}^{m,n}_{p,q}(s)|  = & 
\sum_{i=1}^m \log|G(b_i+\beta_i s;\tau)| + \sum_{i=1}^n \log|G(1+\tau-a_i - 
\alpha_i s;\tau)| \\[1ex]
& -\sum_{i=m+1}^q \log|G(1+\tau-b_i-\beta_i s;\tau)| - 
\sum_{i=n+1}^p \log|G(a_i+\alpha_i s;\tau)| ,
\end{split}
\end{equation}
and, after using the above asymptotic expression for $\log{|G(\sigma+\nu s;\tau)|}$ to treat the terms $\log{|G(b_i+\beta_i s;\tau)|}$ and $\log{|G(a_i+\alpha_i s;\tau)|}$ and 
the corresponding one for $\log{|G(\sigma-\nu s;\tau)|}$ to handle  
$\log{|G(1+\tau-a_i-\alpha_i s;\tau)|}$ and $\log{|G(1+\tau-b_i-\beta_i s;\tau)|}$, we obtain
\begin{equation}
\label{integrand.2}
\log|\mathcal{J}^{m,n}_{p,q}(s)|   =  K_1 \, R^2 \log{R} + K_2 \, R^2 + K_3\, R \log{R}
+ K_4 \, R + K_5 \, \log{R} + \mathcal{O}(1) , 
\end{equation}
where 
\begin{align}
& K_1 = \frac{\cos{2\theta}}{2\tau} \, \Delta_2 , \\[1ex]
& K_2 = \frac{\cos{2\theta}}{2\tau}\left[ \Theta_2 - \left(\frac{3}{2}+\log\tau\right) \, \Delta_2 \right] - \frac{\theta\sin{2\theta}}{2\tau}\, \Delta_2 -\big(1-2[\theta>0]\big)\frac{\pi \sin{2\theta}}{2\tau}\, \Omega_2   , \\[1ex]
& K_3 =   \cos{\theta} \, \Phi_1  - 
\frac{\sin\theta}{\tau} \, \operatorname{Im}\Pi ,\\[1ex]
& K_4 =  \cos\theta\bigg(\Phi_2 
 + \frac{1}{2}\log{2\pi} \, \Delta^\ast\bigg)  
- (1+\log\tau)\left[\cos\theta\, \Phi_1 -\frac{\sin\theta}{\tau}\operatorname{Im}\Pi \right] \\[1ex]
 &  -\theta \left[ \sin\theta \, \Phi_1 + \frac{\cos\theta}{\tau}\operatorname{Im} \Pi \right]  - \frac{\sin\theta}{\tau}\, \operatorname{Im}\Lambda + 
\big(1-2[\theta>0]\big)\pi\left[ \sin\theta\, \Phi_3 
+ \cos\theta \frac{\operatorname{Im}\Upsilon}{\tau} \right], \\[1ex]
& K_5 = \frac{1}{2\tau}\left(\operatorname{Re}D_2 + 2c\operatorname{Re}\Pi + c^2 \Delta_2\right)  - \varkappa \,\left(\operatorname{Re}D_1 + c\Delta_1 \right) + \left(\frac{\tau}{12} + \frac{1}{4} + 
\frac{1}{12\tau}\right) \mathcal{N} ,
\end{align}
and 
\begin{align}
\qquad \quad & \Phi_1 =  \frac{\operatorname{Re}\Pi + c\Delta_2}{\tau} - \varkappa  \, \Delta_1  , \quad 
 \Phi_2 = \frac{\operatorname{Re}\Lambda + c\Theta_2}{\tau}-\varkappa 
\, \Theta_1 ,\quad 
\Phi_3 = \frac{\operatorname{Re}\Upsilon-c\Omega_2}{\tau} + \varkappa \, \Omega_1 , \\[1ex]
& \label{def.N} 
\varkappa =  \frac{1}{2}\left(1+\frac{1}{\tau}\right) , \qquad 
\mathcal{N} = 2(m+n)-(p+q) , \\[1ex]
\label{def.Delta_k}
& \Delta_k = \sum_{i=1}^m \beta_i^k + \sum_{i=1}^n \alpha_i^k - \sum_{i=m+1}^q \beta_i^k - \sum_{i=n+1}^p \alpha_i^k , \quad k = 1,2, \\[1ex] 
& \Theta_k = \sum_{i=1}^m \beta_i^k \log\beta_i + \sum_{i=1}^n \alpha_i^k \log\alpha_i - \sum_{i=m+1}^q \beta_i^k \log\beta_i  - \sum_{i=n+1}^p \alpha_i^k \log\alpha_i , \quad k = 1,2, \\[1ex]
& \Omega_k = \sum_{i=1}^n\alpha_i^k  - \sum_{i=m+1}^q \beta_i^k , \quad k = 1,2, \qquad 
\Omega_k^\dagger = \sum_{i=1}^m\beta_i^k  - \sum_{i=n+1}^p \alpha_i^k , \quad k = 1,2,\\[1ex]
& \Delta^\ast = \sum_{i=1}^q \beta_i - \sum_{i=1}^p \alpha_i ,\quad
\Upsilon =   
 \sum_{i=m+1}^q b_i \beta_i - \sum_{i=1}^n a_i \alpha_i ,\quad 
 \Upsilon^\dagger = \sum_{i=n+1}^p a_i\alpha_i - \sum_{i=1}^m b_i \beta_i , 
\\[1ex]
\label{def.Pi_1}
& \Pi =  \sum_{i=1}^m b_i \beta_i + \sum_{i=1}^n a_i \alpha_i - \sum_{i=m+1}^q b_i \beta_i - \sum_{i=n+1}^p a_i \alpha_i , \\[1ex]
& \Lambda =  \sum_{i=1}^m b_i \beta_i \log\beta_i + \sum_{i=1}^n a_i \alpha_i \log\alpha_i - \sum_{i=m+1}^q b_i \beta_i \log\beta_i - \sum_{i=n+1}^p a_i \alpha_i \log\alpha_i ,\\[1ex]
\label{def.D_k}
& D_k = \sum_{i=1}^m b_i^k + \sum_{i=1}^n a_i^k - \sum_{i=m+1}^q b_i^k - \sum_{i=n+1}^p a_i^k , \quad k = 1,2.
\end{align}

Let us analyse the possibilities according to the dominant asymptotic term when $\protect{|s| = R\to\infty}$.  
\medskip

\noindent \textbf{\sffamily (i)} If $\Delta_2 \neq 0$ the leading term is $R^2\log{R}$, and the integral will converge if $K_1 < 0$.  

\begin{enumerate}
\item If $\Delta_2 > 0$ we  choose a 
contour $\mathcal{C}_{ic\infty}$, which can be deformed \cite{D1,D2} into a contour 
$\mathcal{C}_{(\theta_-,\theta_+)}$ with $-3\pi/4 < \theta_- < -\pi/4$ and 
$\pi/4 < \theta_+ < 3\pi/4$. \smallskip

\item If $\Delta_2 < 0$ we can choose either $\mathcal{C}_{-\infty}$ or 
$\mathcal{C}_{+\infty}$; the contour $\mathcal{C}_{-\infty}$ can be 
deformed into a contour $\mathcal{C}_{(\theta_-,\theta_+)}$ 
with $-\pi < \theta_- < -3\pi/4$ and $3\pi/4 < \theta_+ < \pi$, and the 
contour $\mathcal{C}_{+\infty}$ can be 
deformed into a contour $\mathcal{C}_{(\theta_-,\theta_+)}$ 
with $-\pi/4 < \theta_- < 0$ and $0 < \theta_+ < \pi/4$. 
\end{enumerate}

\medskip

\noindent \textbf{\sffamily (ii)} If $\Delta_2 = 0$ or $\cos{2\theta} = 0$ 
the leading term is $R^2$. If $\Delta_2 = 0$ the coefficient of $R^2$ is 
\begin{equation}
\label{eq.R2}
K_2^\prime = \begin{cases}
\frac{1}{2\tau}\left[  \left(\Theta_2 + 2\pi \operatorname{Re}\varepsilon\right) \cos{2\theta} + \pi   \left( \Omega_2  
- 2\operatorname{Im} \varepsilon \right)   \sin{2\theta} \right], & \quad (\theta > 0), \\[1ex]
\frac{1}{2\tau}\left[  \left(\Theta_2 + 2\pi \operatorname{Re}\varepsilon\right) \cos{2\theta} - \pi  \left( \Omega_2  
+ 2 \operatorname{Im} \varepsilon \right) \sin{2\theta} \right], & \quad (\theta < 0) .
\end{cases}
\end{equation}
The integral will converge for $R \to \infty$ if $K_2^\prime < 0$. 
For the contours $\mathcal{C}_{-\infty}$, $\mathcal{C}_{+\infty}$ and $\mathcal{C}_{ic\infty}$ 
we have $\sin{2\theta} = 0$ for $R \to \infty$.

\begin{enumerate}
\setcounter{enumi}{2}
\item If $\Delta_2 = 0$  and
 $(\Theta_2 + 2\pi \operatorname{Re}\varepsilon) > 0 $ 
we choose the contour $\mathcal{C}_{ic\infty}$; if in addition 
we have ($\mathrm{c_1}$) $\Omega_2 > 0$ and $|\operatorname{Im}\varepsilon| < \Omega_2/2$, the contour
$\mathcal{C}_{ic\infty}$ can be deformed into a contour $\mathcal{C}_{(\theta_-,\theta_+)}$ with
$-3\pi/4 < \theta_- < -\pi/2 $ and $\pi/2 < \theta_+ < 3\pi/4$, while if we have
($\mathrm{c_2}$) $\Omega_2 < 0$ and $|\operatorname{Im}\varepsilon| < -\Omega_2/2$, the contour
$\mathcal{C}_{ic\infty}$ can be deformed into a contour $\mathcal{C}_{(\theta_-,\theta_+)}$ with
$-\pi/2 < \theta_- < -\pi/4 $ and $\pi/4 < \theta_+ < \pi/2$. \smallskip

\item If $\Delta_2 = 0$  and 
$(\Theta_2 + 2\pi \operatorname{Re}\varepsilon) < 0 $
we can choose either $\mathcal{C}_{-\infty}$ or 
$\mathcal{C}_{+\infty}$; if in addition we have ($\mathrm{d_1}$) $\Omega_2 > 0$ and 
$|\operatorname{Im}\varepsilon| < \Omega_2/2$ the contour $\mathcal{C}_{-\infty}$ can 
be deformed into a contour $\mathcal{C}_{(\theta_-,\theta_+)}$ with $-\pi < \theta_- < -3\pi/4$ 
and $3\pi/4 < \theta_+ < \pi$, while if the have ($\mathrm{d_2}$) $\Omega_2 < 0$ 
and $|\operatorname{Im}\varepsilon| < -\Omega_2/2$ the contour $\mathcal{C}_{+\infty}$ can 
be deformed into a contour $\mathcal{C}_{(\theta_-,\theta_+)}$ with $-\pi/4 < \theta_- < 0 $ 
and $0 < \theta_+ < \pi/4$. 
\end{enumerate}

\noindent We have $\cos{2\theta} = 0$ for $R\to \infty$ for 
the contours $\mathcal{C}_{(-\pi/4,\pi/4)}$, $\mathcal{C}_{(-3\pi/4,3\pi/4)}$, $\mathcal{C}_{(-\pi/4,3\pi,4)}$ and $\mathcal{C}_{(-3\pi/4,\pi/4)}$, and the coefficient of $R^2$ in this case is 
\begin{equation}
\label{eq.R2.1}
K_2^\prime = \begin{cases}
\frac{\pi}{2\tau}\left(-\frac{1}{4}\Delta_2 + \Omega_2 - 2\operatorname{Im}\varepsilon\right) , \quad (\theta = \pi/4), \\[1ex] 
\frac{\pi}{2\tau}\left(\frac{3}{4}\Delta_2 - \Omega_2 + 2\operatorname{Im}\varepsilon\right) , \quad (\theta = 3\pi/4), \\[1ex]
\frac{\pi}{2\tau}\left(-\frac{1}{4}\Delta_2 + \Omega_2 + 2\operatorname{Im}\varepsilon\right) , \quad (\theta = -\pi/4), \\[1ex]
\frac{\pi}{2\tau}\left(\frac{3}{4}\Delta_2 - \Omega_2 - 2\operatorname{Im}\varepsilon\right) , \quad (\theta = -3\pi/4) .
\end{cases}
\end{equation}
We will exclude the contours $\mathcal{C}_{(-\pi/4,3\pi/4)}$ and 
$\mathcal{C}_{(-3\pi/4,\pi/4)}$ because they are not the limit of
any deformation like in $(\mathrm{c_1})$, $(\mathrm{c_2})$, $(\mathrm{d_1})$ and $(\mathrm{d_2})$ above. 
Thus, in order to have $K_2^\prime < 0$ when $\cos{2\theta} = 0$, we have 
the following: 
\begin{enumerate}
\setcounter{enumi}{4}
\item If $ \Omega_2 - 3\Delta_2/4 > 0$ and $2|\operatorname{Im}\varepsilon| < \Omega_2-3\Delta_2/4$ we choose $\mathcal{C}_{(-3\pi/4,3\pi/4)}$. \smallskip

\item If $ \Delta_2/4 -\Omega_2 > 0$ and $2|\operatorname{Im}\varepsilon| < \Delta_2/4 -\Omega_2$ we choose 
$\mathcal{C}_{(-\pi/4,\pi/4)}$.   
\end{enumerate}

 \medskip

\noindent \textbf{\sffamily (iii)} If $\Delta_2 = 0$, $\Theta_2 + 2\pi\operatorname{Re}\varepsilon = 0$, 
$\Omega_2 =0$ and $\operatorname{Im}\varepsilon = 0$ the leading 
asymptotic term when $R\to \infty$  is $R\log{R}$ with coefficient 
\begin{equation}
K_3  = \cos\theta \, \Phi_1 
- \frac{\sin\theta}{\tau}  \operatorname{Im}\left(\Pi\right) . 
\end{equation}
Then we have:
\begin{enumerate}
\setcounter{enumi}{6}
\item If $\Phi_1 > 0$ we choose $\mathcal{C}_{-\infty}$, and if in 
addition we have ($\mathrm{g_1}$) $\operatorname{Im}\Pi = 0$ 
the contour $\mathcal{C}_{-\infty}$ 
can be deformed into a contour $\mathcal{C}_{(\theta_-,\theta_+)}$ 
with $-\pi < \theta_- < -\pi/2$ and $\pi/2 < \theta_+ < \pi$. \smallskip

\item If $\Phi_1 < 0$ we choose $\mathcal{C}_{+\infty}$, and if in 
addition we have ($\mathrm{h_1}$) $\operatorname{Im}\Pi = 0$  
the contour $\mathcal{C}_{+\infty}$ 
can be deformed into a contour $\mathcal{C}_{(\theta_-,\theta_+)}$ 
with $-\pi/2 < \theta_- < 0$ and $0 < \theta_+ < \pi/2$. 
\end{enumerate}

\noindent The contour $\mathcal{C}_{ic\infty}$ cannot be used 
because for $\theta^- = -\pi/2$ we 
would need  $\operatorname{Im}\Pi < 0$ and for $\theta^+ = \pi/2$ we 
would need  $\operatorname{Im}\Pi > 0$, and if $\operatorname{Im}\Pi = 0$ we 
have $K_3 = 0$ for $\mathcal{C}_{ic\infty}$. For the contours $\mathcal{C}_{(-\pi/4,\pi/4)}$, 
$\mathcal{C}_{(-3\pi/4,3\pi/4)}$, $\mathcal{C}_{(-\pi/4,3\pi,4)}$ and $\mathcal{C}_{(-3\pi/4,\pi/4)}$, 
we can have $K_2^\prime = 0$ and $K_3 < 0$ only for $\mathcal{C}_{(-\pi/4,\pi/4)}$ and 
$\mathcal{C}_{(-3\pi/4,3\pi/4)}$. 

\begin{enumerate}
\setcounter{enumi}{8}
\item If $\Omega_2 = \Delta_2/4$, $\operatorname{Im}\varepsilon = 0$, 
$\Phi_1 < 0$ and $\tau^{-1}|\operatorname{Im}\Pi | < -\Phi_1$ we choose 
$\mathcal{C}_{(-\pi/4,\pi/4)}$. \smallskip

\item If $\Omega_2 = 3\Delta_2/4$, $\operatorname{Im}\varepsilon = 0$, 
$\Phi_1 > 0$ and $\tau^{-1}|\operatorname{Im}\Pi | < \Phi_1$ we  choose 
$\mathcal{C}_{(-3\pi/4,3\pi/4)}$. 
\end{enumerate}

\medskip

\noindent \textbf{\sffamily (iv)} If $\Delta_2 = 0$, $\Theta_2 + 2\pi\operatorname{Re}\varepsilon = 0$, 
$\Omega_2 =0$, $\operatorname{Im}\varepsilon = 0$, $\Phi_1 = 0$ 
and $\operatorname{Im}\Pi = 0$, the leading term is $R$ with coefficient 
\begin{equation}
K_4^\prime = 
\begin{cases} \left(\Psi_1 - \log|z| - \pi \tau^{-1}\operatorname{Im}\Upsilon\right)\cos\theta + 
\left(\operatorname{arg}z - \tau^{-1}\operatorname{Im}\Lambda - \pi \Phi_3\right) \sin\theta, \; & (\theta > 0) , \\[1ex]
\left(\Psi_1 - \log|z| + \pi \tau^{-1}\operatorname{Im}\Upsilon\right)\cos\theta + 
\left(\operatorname{arg}z - \tau^{-1}\operatorname{Im}\Lambda + \pi \Phi_3\right) \sin\theta, \; & (\theta < 0), 
\end{cases}
\end{equation}
where 
\begin{equation}
\Psi_1 = \tau^{-1} \operatorname{Re}\Lambda - \varkappa \Theta_1 + \frac{1}{2}\log{2\pi}\, \Delta^\ast .
\end{equation}
Thus we have:
\begin{enumerate}
\setcounter{enumi}{10}
\item If  $\left(\Psi_1 - \log|z| \mp \pi \tau^{-1}\operatorname{Im}\Upsilon\right) > 0$, that is, 
$|z| < \operatorname{e}^{\Psi_1- \pi\tau^{-1}|\operatorname{Im}\Upsilon|}$, 
we choose 
$\mathcal{C}_{-\infty}$. \smallskip

\item If $\left(\Psi_1 - \log|z| \mp \pi \tau^{-1}\operatorname{Im}\Upsilon\right) < 0$, that is, $|z| > \operatorname{e}^{\Psi_1 + \pi\tau^{-1}|\operatorname{Im}\Upsilon|}$, we choose 
$\mathcal{C}_{+\infty}$. 
\end{enumerate}

\noindent The above expression for $K_4^\prime$ also 
holds for the contours with $\cos{2\theta} = 0$ for $R \to \infty$ and $\Theta_2 + 2\pi\operatorname{Re}\varepsilon = 0$, $\Omega_2 =0$, $\operatorname{Im}\varepsilon = 0$, $\Phi_1 = 0$ 
and $\operatorname{Im}\Pi = 0$. 
So we have 

\begin{enumerate}
\setcounter{enumi}{12}
\item For the contour $\mathcal{C}_{(-\pi/4,\pi/4)}$ the integral 
converges for 
$$
|z| > \operatorname{e}^{\Psi_1-\pi\Phi_3+|\tau^{-1}\operatorname{Im}\Lambda 
+ \pi\tau^{-1}\operatorname{Im}\Upsilon - \operatorname{arg}z|} . 
$$ 

\item For the 
contour $\mathcal{C}_{(-3\pi/4,3\pi/4)}$ the integral 
converges for 
$$|z| < \operatorname{e}^{\Psi_1+\pi\Phi_3-|\tau^{-1}\operatorname{Im}\Lambda 
- \pi\tau^{-1}\operatorname{Im}\Upsilon - \operatorname{arg}z|} . 
$$
\end{enumerate}

\noindent For the contour $\mathcal{C}_{ic\infty}$ the condition $\Phi_1 = 0$ we have assumed above 
is not needed since
$\cos\theta = 0$ for $R\to \infty$ (so $\cos\theta\, \Phi_1 = 0$ in the expression for $K_3$), and the 
expression for $K_4^\prime$ in this case is 
\begin{equation}
K_4^\prime =   
\begin{cases} -\frac{\pi}{2}\Phi_1 - \tau^{-1} \operatorname{Im}\Lambda - \pi \Phi_3 + 
\operatorname{arg}z , \quad (\theta = \pi/2) ,\\[1ex]
-\frac{\pi}{2}\Phi_1 + \tau^{-1} \operatorname{Im}\Lambda - \pi \Phi_3 - 
\operatorname{arg}z , \quad (\theta = -\pi/2) . 
\end{cases}
\end{equation}
Thus we have

\begin{enumerate}
\setcounter{enumi}{14}
\item If $(\Phi_1 + 2 \Phi_3)\pi/2 > 0$ and $|\operatorname{arg}z - \tau^{-1}\operatorname{Im}\Lambda| 
< (\Phi_1 + 2 \Phi_3)\pi/2$ we  choose $\mathcal{C}_{ic\infty}$. 
\end{enumerate}

\medskip

\noindent \textbf{\sffamily (v)} If $\Delta_2 = 0$, $\Theta_2 + 2\pi\operatorname{Re}\varepsilon = 0$, 
$\Omega_2 =0$, $\operatorname{Im}\varepsilon = 0$, $\Phi_1 = 0$, $\operatorname{Im}\Pi = 0$ 
and $K_4^\prime = 0$, the dominant term is $\log{R}$. 
The convergence requires that $K_5 < -1$, that is,    
\begin{equation}
\label{K5.cond}
\chi + c \Phi_1
<-1 .
\end{equation}
where
\begin{equation}
\chi = \frac{1}{2\tau}\operatorname{Re}D_2  - \varkappa \,\operatorname{Re}D_1 
 + \left(\frac{\tau}{12} +  \frac{1}{4} +
\frac{1}{12\tau}\right) \mathcal{N} . 
\end{equation}
\begin{enumerate}
\setcounter{enumi}{15}
\item For $\mathcal{C}_{-\infty}$ and $\mathcal{C}_{+\infty}$ 
with $|z| = \operatorname{e}^{\Psi_1 - \pi \tau^{-1}|\operatorname{Im}\Upsilon|}$ we have $K_4^\prime = 0$, and since $\Phi_1 = 0$ the inequality above becomes 
$\chi < -1 .$\smallskip

\item For $\mathcal{C}_{ic\infty}$ we don't need $\Phi_1 = 0$ but $\Phi_1 + 2 \Phi_3 = 0$ and 
$\operatorname{arg}z = \tau^{-1} \operatorname{Im}\Lambda$. 
If $\Phi_1 > 0$ the above inequality can be conveniently written as 
$c < \mathcal{Z} $
with 
\begin{equation}
\label{def.Z}
\mathcal{Z} = - \frac{1+\chi}{\Phi_1} , 
\end{equation}
while if $\Phi_1 < 0$ the inequality is $c > \mathcal{Z}$. 
\end{enumerate}

\medskip

\subsection{Discussion and Examples}

At this point it is convenient to introduce some terminology. Let us recall that in the case of the gamma function, the highest-order term when $R \to \infty$ is $R \log{R}$. Therefore, the terms $R^2 \log{R}$ and $R^2$ are characteristic of the double gamma function. We will say that a Fox-Barnes $J$-function is \textit{hyper-unbalanced} if $\Delta_2 \neq 0$ -- the use of the prefix hyper is justified to distinguish this situation from that of the Fox $H$-function. Thus a hyper-unbalanced Fox-Barnes $J$-function exists depending on the choice of a contour $\mathcal{C}_{ic\infty}$ or $\mathcal{C}_{-\infty}/\mathcal{C}_{+\infty}$, as in the examples below.
 We will say that a Fox-Barnes $J$-function is \textit{hyper-balanced} if $\Delta_2 = 0$. Since $\Delta_2 = \Omega_2^\dagger - \Omega_2$, if 
a Fox-Barnes $J$-function is hyper-balanced we have $\Omega_2 =\Omega_2^\dagger $. We say a Fox-Barnes $J$-function is 
 \textit{strongly hyper-balanced} if $\Delta_2 = 0$ and $\Omega_2 = 0$, and \textit{completely hyper-balanced} when $\Delta_2 = 0$, $\Omega_2 = 0$, $\Theta_2 = 0$, and $\varepsilon = 0$. Therefore, for a completely hyper-balanced Fox-Barnes $J$-function, the dominant term when $R\to \infty$ is $R\log{R}$, like with the Fox $H$-functions.
In fact, we will see in the next section that 
Fox $H$-functions are examples
of completely hyper-balanced Fox-Barnes $J$-functions. On the 
other hand, the Kilbas-Saigo function and its Laplace transform 
(see Propositions~\ref{prop.6}~and~\ref{prop.7}) 
are not examples of Fox $H$-functions but are examples
of completely hyper-balanced Fox-Barnes $J$-functions. Note that
the quantities $\Delta_2$, $\Omega_2$ and $\Theta_2$ depend only on 
the scale parameters $\alpha_i$ and $\beta_i$ of the Fox-Barnes $J$-function 
and not on the translation parameters $a_i$ and $b_i$.

Let us see some examples of Fox-Barnes $J$-functions.

\begin{figure}
\begin{center}
\includegraphics[width=11cm]{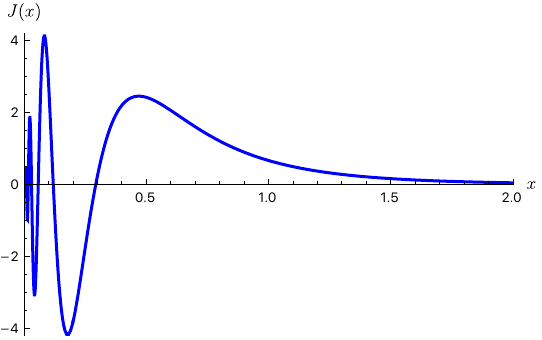}
\caption{Plot of $J(x) = J^{1,0}_{0,1}[x;1|-;(1,2/3)]$ for $x \in (0,2)$. \label{figure_1}}
\end{center}
\end{figure}

\paragraph{Example 1.} We start with the simple example 
\begin{equation}
J^{1,0}_{0,1}\left[x;1\bigg|\begin{matrix} - \\ (1,2/3) \end{matrix}\right] = 
\frac{1}{2\pi i}\int_{\mathcal{C}_1} x^{-s} G(1+2s/3)\, ds .
\end{equation} 
In this case $\tau = 1$ and $\varepsilon = 0$. We have $\Delta_2 = (2/3)^2 > 0$ and so we choose a contour $\mathcal{C}_1 = \mathcal{C}_{ic\infty} = (c-i\infty,c+i\infty)$, and since $G(1+2s/3)$ is an analytic function, we can choose $c$ as any real number. In Figure~\ref{figure_1} we have a plot of $J(x) = J^{1,0}_{(0,1)}[x;1|-;(1,2/3)]$. 
We observe that  $J(x) = J^{1,0}_{0,1}[x;1|-;(1,2/3)]$ is indeed a Mellin-Barnes $K$-function \cite{Kusnetsov}, and this will always happen when we have only one double gamma function in the integrand, that is, a Fox-Barnes $J$-function involving only one double gamma function reduces to a Mellin-Barnes $K$-function  after a change of the integration variable of the form $s \to k s$ (see property (iii) in Proposition~\ref{prop.2}), in this case $k=3/2$. The plot in Figure~\ref{figure_1} was created using Mathematica 14.2 with numerical contour integration and the native function \texttt{BarnesG[z]} for the evaluation 
of the Barnes $G$-function. 

\begin{figure}
\begin{center}
\includegraphics[width=11cm]{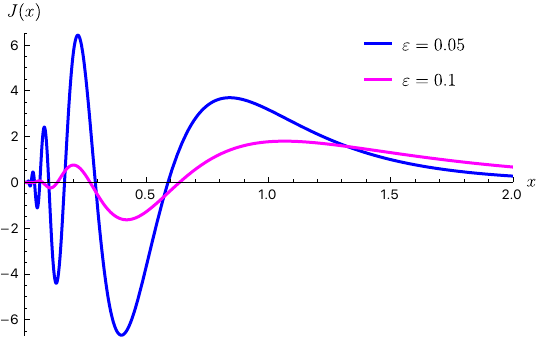}
\caption{Plot of $J(x) = J^{1,0}_{0,1}[x;2,\varepsilon|-;(1,2/3)]$ with 
$\varepsilon = 0.05$ (blue) and $\varepsilon = 0.1$ (magenta) for $x \in (0,2)$. \label{figure_1_bis}}
\end{center}
\end{figure}

\paragraph{Example 2.} An interesting variation of the previous example is 
\begin{equation}
J^{1,0}_{0,1}\left[x;2,\varepsilon \bigg|\begin{matrix} - \\ (1,2/3) \end{matrix}\right] = 
\frac{1}{2\pi i}\int_{\mathcal{C}_1} x^{-s} \operatorname{e}^{\pi \varepsilon s^2/2}G(1+2s/3;2) \, ds 
\end{equation}
with $\varepsilon \neq 0$ and $\tau = 2$. 
Note that $\Delta_2$ has not changed, so we can still use the same contour of Example~1. In Figure~\ref{figure_1_bis} we have a plot of this new situation
for $\varepsilon = 0.05$ (plot in blue) and $\varepsilon = 0.1$ (plot in magenta). The evaluation of the 
values of $G(1+2s/3;2)$ was done using Mathematica 14.2 and numerical integration of the expression in eq.\eqref{double.gamma.int.rep}. Comparing 
Figures~\ref{figure_1}~and~\ref{figure_1_bis}, no clear interpretation 
of the role of the parameter $\varepsilon$ is seen beyond its purely technical role (see also Remark~2).

\begin{figure}
\begin{center}
\includegraphics[width=11cm]{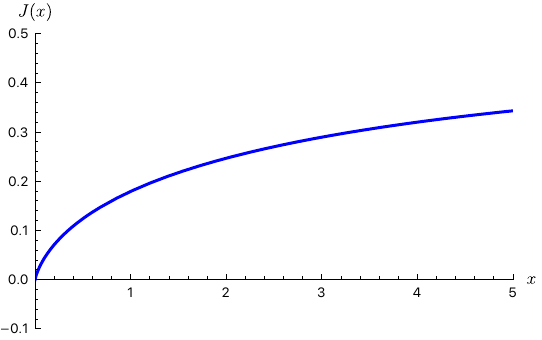}
\caption{Plot of $J(x) = J^{1,1}_{2,1}[x,1|-;(1,3),(1,1);(1,1/2)]$ for $x \in (0,5)$. \label{figure_2}}
\end{center}
\end{figure}

\paragraph{Example 3.} Another example is 
\begin{equation}
J^{1,1}_{2,1}\left[x,1\bigg|\begin{matrix} (1,3), (1,1) \\ (1,1/2) \end{matrix}\right] = 
\frac{1}{2\pi i}\int_{\mathcal{C}_2} x^{-s} \frac{G(1+s/2)G(1-3s)}{G(1+s)} \, ds .
\end{equation}
In this case we also have $\tau = 1$ and $\varepsilon = 0$. We have $\Delta_2 = (1/2)^2 + 3^2 - 1^2  > 0$, so we choose a contour $\mathcal{C}_{ic\infty}$. The integrand has poles at $\grave{\mys}_{m,n} = -(m+n+1)$, or $\grave{\mys}_k = -(k+1)$ for $k=0,1,2,\ldots$ of order $k+1$. So we choose $c$ such that $c > -1$. 
In Figure~\ref{figure_2} we have a plot of $J(x) = J^{1,1}_{(2,1)}[x,1|-;(1,3),(1,1);(1,1/2)]$, also created using Mathematica 14.2. 

\paragraph{Example 4.} A simple but interesting example is 
\begin{equation}
J^{0,0}_{1,0}\left[x;1\bigg|\begin{matrix} (0,1/\nu) \\ - \end{matrix}\right] = 
\frac{1}{2\pi i}\int_{\mathcal{C}_3} \frac{x^{-s}}{G(s/\nu)} \, ds = 
\frac{\nu}{2\pi i}\int_{\mathcal{C}_3}\frac{(x^\nu)^{-s}}{G(s)}\, ds = 
\nu J^{0,0}_{1,0}\left[x^\nu;1\bigg|\begin{matrix} (0,1) \\ - \end{matrix}\right].
\end{equation}
In this case $\Delta_2 = -1/\nu^2 < 0$ and so we can choose $\mathcal{C}_{-\infty}$ 
or $\mathcal{C}_{+\infty}$. The poles of the integrand comes from the zeros of $G(s)$ in the denominator, that is, $\grave{s}_{m,n} = -(m+n)$ for $m,n = 0,1,2,\ldots$, or $\grave{s}_k = -k$ of order $k+1$ for $k=0,1,2,\ldots$. Thus the choice giving a non-null function is $\mathcal{C}_{-\infty}$. Calculating the values of the residues at $\grave{s}_k$ 
becomes increasingly intricate as the order of the poles increases, 
but it is not our intention to perform this calculation explicitly, only to understand the form of the resulting series. To do this, we recall that at 
a neighbourhood of $\grave{s}_k$ we can write $G(s) = (s+k)^{k+1}H_k(s)$, 
and defining $s+k = t$, we can write the integrand as $x^k \varphi_k(t,x)/t^{k+1}$ 
with $\varphi_k(t,x) = \operatorname{e}^{-t\log{x}}/H_k(t-k)$. Thus the residue at $\grave{s}_k = -k$ is 
\begin{equation}
\underset{s=-k}{\operatorname{Res}}\frac{x^{-s}}{G(s)} = x^k 
\underset{t=0}{\operatorname{Res}}\frac{\varphi_k(t,x)}{t^{k+1}} = 
\frac{x^k}{k!}\frac{d^k}{dt^k}\left[\frac{\operatorname{e}^{-t\log{x}}}{H_k(t-k)}\right]_{t=0} = x^k \sum_{n=0}^k c_{k,n} (\log{x})^k ,
\end{equation}
for some constants $c_{k,n}$. Reintroducing the parameter $\nu$ and using the residue theorem we obtain that $J^{0,0}_{1,0}[x;1|(0,1/\nu);-]$ is of the 
form 
\begin{equation}
J^{0,0}_{1,0}\left[x;1\bigg|\begin{matrix} (0,1/\nu) \\ - \end{matrix}\right] =
\sum_{k=0}^\infty \sum_{n=0}^k \nu^{k+1} c_{k,n} x^{\nu n} \log^k{x} .
\end{equation}
Thus $J^{0,0}_{1,0}[x;1|(0,1/\nu);-]$ can be represented by a power-logarithmic series.

\medskip

We also recall that different contours define different functions if they cannot be continuously deformed into each other. Therefore, because of its relationship with the Mellin transform, we will take 
the contour $\mathcal{C}_{ic\infty}$ as the primary choice for the definition of the Fox-Barnes $J$-function. It is convenient to collect the conditions seen above that lead to the contour $\mathcal{C}_{ic\infty}$ with the help of a denomination:

\enlargethispage{2ex}

\begin{defin} \label{def.1} A vertical line (VL) condition for the Fox-Barnes $J$-function is any of the following conditions: 
\begin{enumerate}
\item $\Delta_2 > 0$.
\item $\Delta_2 = 0$ and $\Theta_2 + 2\pi \operatorname{Re}\varepsilon > 0$.
\item $\Delta_2 =0$, $\Theta_2 + 2 \pi \operatorname{Re}\varepsilon =0$, $\Omega_2 = 0$, $\operatorname{Im}\varepsilon = 0$, $\operatorname{Im}\Pi = 0$, $\Phi_1 + 2 \Phi_3 > 0$ and \\
$|\operatorname{arg}z - \tau^{-1}\operatorname{Im}\Lambda| < (\Phi_1 + 2 \Phi_3)\pi/2$.
\item $\Delta_2 =0$, $\Theta_2 + 2 \pi \operatorname{Re}\varepsilon =0$, $\Omega_2 = 0$, $\operatorname{Im}\varepsilon = 0$, $\operatorname{Im}\Pi = 0$, $\Phi_1 + 2 \Phi_3 =0$, $\operatorname{arg}z = \tau^{-1}\operatorname{Im}\Lambda$ and $\chi + c \Phi_1 < -1$.  
\end{enumerate}
\end{defin}

\section{A particular case: the Fox $H$-function}
\label{sec.4}

Let us recall the definition of Fox $H$-function, that is, 
\begin{equation}
\label{alternative.fox.H}
H^{m,n}_{p,q}\left[ z\bigg| \begin{matrix} (
{a}_i,{\alpha}_i)_{1,p} \\ ({b}_i,\beta_i)_{1,q} \end{matrix}\right] = 
\frac{1}{2\pi i} \int_{\gamma} \mathcal{H}^{m,n}_{p,q}(s) z^{-s}\, ds ,
\end{equation}
with $\mathcal{H}^{m,n}_{p,q}(s)$ given by 
\begin{equation}
\label{alternative.fox.H.integrand}
\mathcal{H}^{m,n}_{p,q}(s) = \frac{\prod_{i=1}^{m} \Gamma({b}_i+{\beta}_i s) \cdot \prod_{i=1}^{n} \Gamma(1-{a}_i-{\alpha}_i s)}{ 
\prod_{i=m+1}^{q} \Gamma(1-{b}_i-{\beta}_i s)\cdot \prod_{i=n+1}^{p} \Gamma({a}_i + {\alpha}_i s) } , 
\end{equation}
and $\gamma$ is one of the contours $\mathcal{C}_{-\infty}$, $\mathcal{C}_{+\infty}$ and 
$\mathcal{C}_{ic\infty}$ depending on the parameters of the $H$-function. 

We want to show the following 

\bigskip

\begin{prop} 
\label{prop.1}
Consider the Fox-Barnes $J$-function as in eq.\eqref{fox-barnes} and eq.\eqref{fox-barnes.integrand} with $\varepsilon = 0$. Suppose that 
\begin{equation} 
\label{cond.1}
p = q = m + n 
\end{equation} 
and  
\begin{equation}
\label{cond.2}
\begin{aligned}
& \alpha_i = \beta_{m+i} , \qquad i = 1, \ldots, n , \\[1ex]
& \beta_i = \alpha_{n+i} , \qquad i = 1, \ldots, m .
\end{aligned}
\end{equation}
Suppose also that there exists $m_0$ and $n_0$ with $1 \leq m_0 \leq m$ and $1 \leq n_0 \leq n$ such that
\begin{equation}
\label{cond.3}
\begin{aligned}
& a_i = \begin{cases} b_{m+i} - 1 , \qquad & i = 1,\ldots, n_0 , \\[1ex]
b_{m+i} + 1, \qquad & i = n_0 + 1, \ldots, n , \end{cases} \\[1ex]
& b_i = \begin{cases} a_{n+i}+1 , \qquad & i = 1,\ldots, m_0 , \\[1ex]
 a_{n+i} - 1 , \qquad & i = m_0 + 1, \ldots, m . \end{cases} 
\end{aligned}
\end{equation}
Then we have
\begin{equation}
\label{L.2.Fox}
 J^{m,n}_{p,q}\left[ z;\tau\, \bigg| \, \begin{matrix} (a_i,\alpha_i)_{1,p} \\ (b_i,\beta_i)_{1,q}\end{matrix}\right]  
=   
\tau H^{m_0,n_0}_{p_0,q_0}\left[z^\tau \bigg| \begin{matrix} (
\tilde{a}_i,\tilde{\alpha}_i)_{1,p_0}  \\ (\tilde{b}_i,\tilde{\beta}_i)_{1,q_0}\end{matrix}\right] 
\end{equation}
where 
\begin{align}
\label{def.p0.q0}
& p_0 = n_0 + m - m_0 , \qquad q_0 = m_0 + n - n_0 , \\[1ex] 
\label{def.tilde.a}
& \tilde{a}_k = \begin{cases} a_k/\tau , \; & k = 1,\ldots, n_0 \\[1ex]
 \left(a_{n-n_0+m_0+k}-1\right)/\tau , \;& k = n_0+1,\ldots,  p_0 \end{cases} \\[1ex] 
 \label{def.tilde.b}
& \tilde{b}_k = \begin{cases} \left(b_k -1\right)/\tau , \; & k = 1,\ldots, m_0 \\[1ex]
 b_{m-m_0+n_0+k}/\tau , \;& k = m_0+1,\ldots, q_0 \end{cases} \\[1ex] 
 \label{def.tilde.alpha}
& \tilde{\alpha}_k = \begin{cases} \alpha_k   , \; & k = 1,\ldots,n_0 \\[1ex]
\alpha_{n-n_0 + m_0 + k} , \; & k = n_0 + 1,\ldots, p_0 \end{cases} \\[1ex]
\label{def.tilde.beta}
& \tilde{\beta}_k = \begin{cases} \beta_k   , \; & k = 1,\ldots,m_0 \\[1ex]
\beta_{m-m_0 + n_0 + k}   , \; & k = m_0 + 1,\ldots, q_0 \end{cases} 
\end{align}
\end{prop}

 \bigskip
 
\begin{proof} Let us  write each term in the definition of $\mathcal{J}^{m,n}_{p,q}(s)$ in eq.\eqref{fox-barnes.integrand} as  
\begin{align}
\label{proof.1.1}
& \prod_{i=1}^m G(b_i+\beta_i s;\tau) = 
\prod_{i=1}^{m_0}G(b_i+\beta_i s;\tau)\prod_{i=m_0+1}^m G(b_i+\beta_i s;\tau) , \\[1ex]
\label{proof.1.2}
& \prod_{i=1}^n G(1+\tau -a_i-\alpha_i s;\tau) =
\prod_{i=1}^{n_0}G(1+\tau-a_i-\alpha_i s;\tau)\!\!\! \prod_{i=n_0+1}^n G(1+\tau-a_i-\alpha_i s;\tau) , \\[1ex]
& \prod_{i=m+1}^q G(1+\tau -b_i-\beta_i s;\tau) = 
\prod_{i=m+1}^{m+n_0}G(1+\tau-b_i-\beta_i s;\tau) \!\!\!\!  \prod_{i=m+n_0+1}^q G(1+\tau-b_i-\beta_i s;\tau)   \\[1ex]
\label{proof.1.3}
& \phantom{G(1+\tau -b_i-\beta_i s;\tau)} = 
\prod_{i=1}^{n_0}G(1+\tau-b_{m+i}-\beta_{m+i} s;\tau) \!\!\!\! \prod_{i=n_0+1}^n G(1+\tau-b_{m+i}-\beta_{m+i} s;\tau) \\[1ex]
& \prod_{i=n+1}^p G(a_i +\alpha_i s;\tau) = 
\prod_{i=n+1}^{n+m_0}G(a_i+\alpha_i s;\tau)\!\!\! \prod_{i=n+m_0+1}^p G(a_i+\alpha_i s;\tau) \\[1ex]
\label{proof.1.4}
& \phantom{\prod_{j=n+1}^p G(a_j +\alpha_j s;\tau)} = 
\prod_{i=1}^{m_0}G(a_{n+i}+\alpha_{n+i} s;\tau)\prod_{i=m_0+1}^m G(a_{n+i}+\alpha_{n+i} s;\tau) .
\end{align}
where in \eqref{proof.1.3} and \eqref{proof.1.4} we have also used eq.\eqref{cond.1}. 
Using also \eqref{cond.2} and \eqref{cond.3} we can rewrite \eqref{proof.1.1}, \eqref{proof.1.2}, \eqref{proof.1.3} and \eqref{proof.1.4} as 
 \begin{align}
\label{proof.1.1.1}
& \prod_{i=1}^m G(b_i+\beta_i s;\tau) = 
\prod_{i=1}^{m_0}G(b_i+\beta_i s;\tau)\prod_{i=m_0+1}^m G(a_{n+i}-1+\alpha_{n+i} s;\tau) , \\[1ex]
\label{proof.1.2.1}
& \prod_{i=1}^n G(1+\tau -a_i-\alpha_i s;\tau) =
\prod_{i=1}^{n_0}G(1+\tau-a_i-\alpha_i s;\tau)\prod_{i=n_0+1}^n G(\tau-b_{m+i}-\beta_{m+i} s;\tau) , \\[1ex]
\label{proof.1.3.1}
& \prod_{i=m+1}^q G(1+\tau -b_i-\beta_i s;\tau) = 
\prod_{i=1}^{n_0}G(\tau-a_i-\alpha_i s;\tau)\prod_{i=n_0+1}^n G(1+\tau-b_{m+i}-\beta_{m+i} s;\tau) \\[1ex]
\label{proof.1.4.1}
& \prod_{i=n+1}^p G(a_i +\alpha_i s;\tau) = 
\prod_{i=1}^{m_0}G(b_i-1+\beta_i s;\tau)\prod_{i=m_0+1}^m G(a_{n+i}+\alpha_{n+i} s;\tau) .
\end{align}
Next we use the above expressions in eq.\eqref{fox-barnes.integrand} and
regroup the terms according to the same product notation, that is, 
\begin{equation}
\begin{split}
\mathcal{J}^{m,n}_{p,q}(s) = &
\prod_{i=1}^{m_0}\frac{G(b_i+\beta_i s;\tau)}{G(b_i-1+\beta_i s;\tau)} 
\cdot \prod_{i=1}^{n_0} \frac{G(1+\tau-a_i-\alpha_i s;\tau)}{G(\tau-a_j -\alpha_j s;\tau)} \\[1ex]
&\cdot \prod_{i=m_0+1}^m\frac{G(a_{n+1}-1+\alpha_{n+i}s;\tau)}{G(a_{n+i}+\alpha_{n+i}s;\tau)} 
\cdot \prod_{i=n_0+1}^n \frac{G(\tau-b_{m+i}-\beta_{m+i}s;\tau)}{G(1+\tau-b_{m+i}-\beta_{m+i}s;\tau)} .
\end{split}
\end{equation}
Using eq.\eqref{ap.B.general.G.tau} we rewrite this as 
\begin{equation}
\mathcal{J}^{m,n}_{p,q}(s) = \frac{\displaystyle \prod_{i=1}^{m_0} \Gamma\left(\frac{b_i-1+\beta_i s}{\tau}\right) 
\cdot \prod_{i=1}^{n_0} \Gamma\left(\frac{\tau-a_i-\alpha_i s}{\tau}\right)}{ \displaystyle 
\prod_{i=m_0+1}^m \Gamma\left(\frac{a_{n+i}-1+\alpha_{n+i}s}{\tau}\right) \cdot 
\prod_{i=n_0+1}^n \Gamma\left(\frac{\tau-b_{m+i}-\beta_{m+i}s}{\tau}\right) } . 
\end{equation}
Defining $\tilde{a}_k$, $\tilde{b}_k$, $\tilde{\alpha}_k$ and $\tilde{\beta}_k$ as in 
\eqref{def.tilde.a}, \eqref{def.tilde.b}, \eqref{def.tilde.alpha} and \eqref{def.tilde.beta}, respectively, 
and $p_0$ and $q_0$ as in \eqref{def.p0.q0}, we rewrite the above expression as 
\begin{equation}
\mathcal{J}^{m,n}_{p,q}(s) = 
\frac{\prod_{i=1}^{m_0}\Gamma(\tilde{b}_i + \tilde{\beta}_i s/\tau)\cdot \prod_{i=1}^{n_0}(1-\tilde{a}_i -\tilde{\alpha}_i s/\tau)}{\prod_{i=m_0+1}^{q_0} \Gamma(1-\tilde{b}_i - \tilde{\beta}_i s/\tau) \cdot 
\prod_{i=n_0+1}^{p_0} \Gamma(\tilde{a}_i + \tilde{\alpha}_j s/\tau)}  . 
\end{equation}
Finally, making the change of variable $s \to \tau s$ and 
recalling eq.\eqref{alternative.fox.H.integrand}, we recognize 
the definition of the Fox $H$-function in eq.\eqref{L.2.Fox}. 
\end{proof}

It is interesting to note that a Fox-Barnes $J$-function satisfying the conditions \eqref{cond.1} and \eqref{cond.2} is completely hyper-balanced. Thus a Fox $H$-function has the particularity of being a completely hyper-balanced function with the additional conditions \eqref{cond.3} on the parameters 
$a_i$ ($i=1,\ldots,n$) and $b_i$ ($i=1,\ldots,m$). 

\subsection{Existence}

When the Fox-Barnes $J$-function reduces to the Fox $H$-function 
as in Proposition~\ref{prop.1}, 
it is straightforward to see that if the
condition in eq.\eqref{cond.2} is satisfied, we have
\begin{equation}
\label{eq.auxiliar.1}
\begin{aligned}
& \Delta_2 = \Delta_1 = 0 , \qquad 
  \Theta_2 = \Theta_1 = 0 , \\[1ex]
& \Omega_2 = \Omega_1 = 0 , \qquad 
  \Delta^\ast = 0 , \qquad 
 \mathcal{N} = 0 , 
\end{aligned}
\end{equation}
and using eq.\eqref{cond.3},  
\begin{equation}
\label{eq.auxiliar.2}
\begin{aligned}
& \Pi = \tilde{\Delta}^\ast , \qquad 
  \Lambda = \log\tilde{\mathcal{H}} , \qquad 
 \Upsilon = \tilde{\mathcal{Y}} , \\[1ex]
& D_1 = (q_0 - p_0) , \qquad 
  D_2 = (q_0 - p_0) + 2 \tau \tilde{\delta}^\ast ,
\end{aligned}
\end{equation}
where 
\begin{equation}
\label{eq.auxiliar.3}
\begin{aligned}
& \tilde{\Delta}^\ast  = \sum_{i=1}^{q_0} \tilde{\beta}_i - \sum_{i=1}^{p_0} \tilde{\alpha}_i , \qquad 
  \tilde{\mathcal{H}} = \prod_{i=1}^{q_0} \tilde{\beta}_i^{\tilde{\beta}_i} 
\prod_{i=1}^{p_0} \tilde{\alpha}_i^{-\tilde{\alpha}_i} , \\[1ex]
& \tilde{\mathcal{Y}} = \sum_{i=1}^{n_0} \tilde{\alpha}_i 
- \sum_{i=m_0+1}^{q_0}\tilde{\beta}_i , \qquad 
  \tilde{\delta}^\ast  = \sum_{i=1}^{q_0} \tilde{b}_i - \sum_{i=1}^{p_0} \tilde{a}_i . 
\end{aligned}
\end{equation}
Consequently we also have 
\begin{equation}
\begin{aligned}
& \operatorname{Im}\Pi = 0 , \qquad 
  \operatorname{Im}\Lambda = 0 , \qquad 
  \operatorname{Im}\Upsilon = 0 .
\end{aligned}
\end{equation}

Using the above expression in eq.\eqref{integrand.2} we obtain 
\begin{align}
& K_1 = 0 , \qquad 
  K_2 = 0 , \qquad 
  K_3 = \tau^{-1} \cos\theta \, \tilde{\Delta}^\ast , \\[1ex]
& K_4 = \tau^{-1} \cos\theta\left[ \log\tilde{\mathcal{H}} - 
(1+\log\tau)\tilde{\Delta}^\ast \right] 
 - \tau^{-1}\sin\theta\left[ \theta\tilde{\Delta}^\ast 
-  (1-2[\theta>0])\pi  \tilde{\mathcal{Y}} \right] , \\[1ex]
& K_5 = \frac{1}{2}(p_0-q_0) + \operatorname{Re} \tilde{\delta}^\ast  + \frac{c}{\tau}\, \tilde{\Delta}^\ast . 
\end{align}

Let us analyse the possibilities according to the dominant asymptotic term when $|s| = R \to \infty$ 
for this particular case. 

\medskip

\noindent \textbf{\sffamily (i)} Since $K_1 = K_2 = 0$, and $\varepsilon = 0$, which gives $K_2^\prime = 0$, the 
dominant asymptotic term when $|s| = R \to \infty$ is $R\log{R}$ with 
coefficient $K_3$. The integral will converge if $K_3 < 0$.  Thus
\begin{enumerate}
\setcounter{enumi}{1}
\item If $\tilde{\Delta}^\ast  > 0$ we  choose the 
contour $\mathcal{C}_{-\infty}$. \smallskip

\item If $\tilde{\Delta}^\ast < 0$ we choose the contour 
$\mathcal{C}_{+\infty}$. 
\end{enumerate}

\medskip

\noindent \textbf{\sffamily (ii)} If $\tilde{\Delta}^\ast = 0$ or $\cos{\theta} = 0$ 
the leading term is $R$ with 
$$
K_4^\prime = K_4 + \tau(-\cos\theta \,\log{|z|} + \sin\theta \, \operatorname{arg}z),
$$ 
where we recall that in this case we changed  $z \to z^\tau$ in eq.\eqref{integrand}. 
Thus 
\begin{enumerate}
\setcounter{enumi}{2}
\item For $\tilde{\Delta}^\ast  = 0$ and $|z| < \big(\tilde{\mathcal{H}}\big)^{\tau^{-2}}$ 
we choose $\mathcal{C}_{-\infty}$. \smallskip

\item For $\tilde{\Delta}^\ast = 0$ and $|z| > \big(\tilde{\mathcal{H}}\big)^{\tau^{-2}}$ we choose $\mathcal{C}_{+\infty}$. \smallskip

\item For $\tilde{\Delta}^\ast  = 0$, 
$\tilde{\mathcal{Y}} > 0$ and $|\operatorname{arg}z| <  \tau^{-2}\tilde{\mathcal{Y}}\, \pi$  
we choose $\mathcal{C}_{i c \infty}$. 
\end{enumerate}

\noindent However, for $\mathcal{C}_{ic\infty}$ we don't need to suppose 
$\tilde{\Delta}^\ast = 0$ since in this case $\cos{\theta} = 0$. With $\cos\theta = 0$ the condition 
$K_4^\prime < 0$ implies that 

\begin{enumerate}
\setcounter{enumi}{5}
\item for $|\operatorname{arg}z| <  \tau^{-2} \tilde{\mathcal{A}}\, \pi/2$ 
we choose $\mathcal{C}_{i c\infty}$, where 
\begin{equation}
\tilde{\mathcal{A}} = \tilde{\Delta}^\ast + 2 \tilde{\mathcal{Y}} = 
\sum_{i=1}^{m_0} \tilde{\beta}_i + \sum_{i=1}^{n_0}\tilde{\alpha}_i 
- \sum_{i=m_0 + 1}^{q_0}\tilde{\beta}_i - \sum_{i=n_0+1}^{p_0} \tilde{\alpha}_i . 
\end{equation}
\end{enumerate}

\medskip

\noindent \textbf{\sffamily (iii)} If all the previous terms vanish, the dominant term is $\log{R}$, 
and convergence requires that $K_5 < -1$. Thus 
\begin{enumerate}
\setcounter{enumi}{6}
\item We can choose $\mathcal{C}_{ic\infty}$ if 
$|\operatorname{arg}z| = \tau^{-2} \tilde{\mathcal{A}}\, \pi/2$ and 
\begin{equation}
 \frac{1}{2}(p_0-q_0) + \operatorname{Re} \tilde{\delta}^\ast  + \frac{c}{\tau}\, \tilde{\Delta}^\ast < -1. 
\end{equation}

\item We can choose $\mathcal{C}_{-\infty}$ or $\mathcal{C}_{+\infty}$ if 
$\tilde{\Delta}^\ast  = 0$, $|z| =  \big(\tilde{\mathcal{H}}\big)^{\tau^{-2}}$ and 
\begin{equation}
 \frac{1}{2}(p_0-q_0) + \operatorname{Re} \tilde{\delta}^\ast   < -1. 
\end{equation}
\end{enumerate}

\subsection{Another possible particular case}

In the proof of Proposition~\ref{prop.1} we have used the functional relation eq.\eqref{ap.B.general.G.tau} 
since the parameters $a_i$ and $b_i$ are supposed to be related by eq.\eqref{cond.3}. On the other hand, 
since in the case of the double gamma function  we have a second functional relation eq.\eqref{ap.B.general.G.tau.2}, 
we can expect that by replacing $1$ with $\tau$ in eq.\eqref{cond.3} we can obtain another 
relation between the Fox-Barnes $J$-function and the Fox $H$-function.

\bigskip

\begin{prop}.
\label{prop.1.5}
 Consider the Fox-Barnes $J$-function as in eq.\eqref{fox-barnes} and eq.\eqref{fox-barnes.integrand} with $\varepsilon = 0$. Suppose that 
\begin{equation} 
p = q = m + n 
\end{equation} 
and  
\begin{equation}
\begin{aligned}
& \alpha_i = \beta_{m+i} , \qquad i = 1, \ldots, n , \\[1ex]
& \beta_i = \alpha_{n+i} , \qquad i = 1, \ldots, m .
\end{aligned}
\end{equation}
Suppose also that there exists $m_0$ and $n_0$ with $1 \leq m_0 \leq m$ and $1 \leq n_0 \leq n$ such that
\begin{equation}
\label{cond.3.1}
\begin{aligned}
& a_i = \begin{cases} b_{m+i} - \tau , \qquad & i = 1,\ldots, n_0 , \\[1ex]
b_{m+i} + \tau, \qquad & i = n_0 + 1, \ldots, n , \end{cases} \\[1ex]
& b_i = \begin{cases} a_{n+i}+\tau , \qquad & i = 1,\ldots, m_0 , \\[1ex]
 a_{n+i} - \tau , \qquad & i = m_0 + 1, \ldots, m . \end{cases} 
\end{aligned}
\end{equation}
Then we have
\begin{equation}
\label{L.2.1.Fox}
 J^{m,n}_{p,q}\left[ z;\tau\, \bigg| \, \begin{matrix} (a_i,\alpha_i)_{1,p} \\ (b_i,\beta_i)_{1,q} \end{matrix}\right]  
 = 2\pi^{(\tau-1)\mathcal{N}_0/2} \tau^{\nu_0 - \bar{\delta}^\ast}  
H^{m_0,n_0}_{p_0,q_0}\left[ \tau^{\bar{\Delta}^\ast} \, z\bigg| \begin{matrix} (
\bar{a}_i,\bar{\alpha}_i)_{1,p_0}   \\ (\bar{b}_i,\bar{\beta}_i)_{1,q_0}  \end{matrix}\right] 
\end{equation}
where 
\begin{align}
& p_0 = n_0 + m - m_0 , \qquad q_0 = m_0 + n - n_0 , \\[1ex] 
& \mathcal{N}_0 = 2(m_0 + n_0) - (p_0+q_0)  , \qquad \nu_0 = \frac{q_0 - p_0}{2} , \\[1ex]
& \bar{\Delta}^\ast = \sum_{k=1}^{q_0} \bar{\beta}_k - \sum_{k=1}^{p_0} \bar{\alpha}_k , \qquad 
\bar{\delta}^\ast = \sum_{k=1}^{q_0} \bar{b}_k - \sum_{k=1}^{p_0} \bar{a}_k , \\[1ex]
& \bar{a}_k = \begin{cases} a_k , \; & k = 1,\ldots, n_0 \\[1ex]
a_{n-n_0+m_0+k}-\tau , ;& k = n_0+1,\ldots, p_0 \end{cases} \\[1ex] 
& \bar{b}_k = \begin{cases} b_k -\tau , \; & k = 1,\ldots, m_0 \\[1ex]
b_{m-m_0+n_0+k}  , ;& k = m_0+1,\ldots, q_0 \end{cases} \\[1ex] 
& \bar{\alpha}_k = \begin{cases} \alpha_k , \; & k = 1,\ldots,n_0 \\[1ex]
\alpha_{n-n_0 + m_0 + k} , \; & k = n_0 + 1,\ldots, p_0 \end{cases} \\[1ex]
& \bar{\beta}_k = \begin{cases} \beta_k , \; & k = 1,\ldots,m_0 \\[1ex]
\beta_{m-m_0 + n_0 + k} , \; & k = m_0 + 1,\ldots, q_0 \end{cases} 
\end{align}
\end{prop}

 \bigskip
 
\begin{proof} The proof is analogous to the previous proposition using 
eq.\eqref{cond.3.1} instead of eq.\eqref{cond.3}, eq.\eqref{ap.B.general.G.tau.2} 
instead of eq.\eqref{ap.B.general.G.tau}, the definitions of 
$\bar{a}_k$, $\bar{b}_k$, $\bar{\alpha}_k$ and $\bar{\beta}_k$ instead 
of $\tilde{a}_k$, $\tilde{b}_k$, $\tilde{\alpha}_k$ and $\tilde{\beta}_k$, 
and the quantities $p_0$, $q_0$, $c_0^\ast$, $\nu_0$, $\bar{\Delta}$ 
and $\bar{\delta}$ as in Proposition~\ref{prop.1.5} 
\end{proof}

\bigskip

\paragraph{Remark 1.} The results in eq.\eqref{eq.auxiliar.1} still hold in the case in Proposition~\ref{prop.1.5}. 
However, eq.\eqref{eq.auxiliar.2} have to be replaced by 
\begin{equation}
\label{eq.auxiliar.2.1}
\begin{aligned}
& \Pi = \tau \bar{\Delta}^\ast , \qquad 
  \Lambda = \tau \log\bar{\mathcal{H}} , \qquad 
  \Upsilon = \tau \bar{\mathcal{Y}} , \\[1ex]
& D_1 = \tau (q_0 - p_0) , \qquad 
  D_2 = \tau^2 (q_0 - p_0) + 2 \tau \bar{\delta}^\ast ,
\end{aligned}
\end{equation}
where $\bar{\Delta}^\ast$, $\bar{\mathcal{H}}$, $\bar{\mathcal{Y}}$ and 
$\bar{\delta}$ are defined as in eq.\eqref{eq.auxiliar.3} but with $\tilde{\alpha}_i$, $\tilde{\beta}_i$, 
$\tilde{a}_i$ and $\tilde{b}_i$ replaced by $\bar{\alpha}_i$, $\bar{\beta}_i$, $\bar{a}_i$ and $\bar{b}_i$.

\section{Properties}
\label{sec.5}

The Fox-Barnes $J$-function has some properties completely analogous to the 
Fox $H$-function, that is: 
\begin{description}
\item[(i)] It is symmetric in the set of pairs 
$\{(a_1,\alpha_1),\ldots,(a_n,\alpha_n)\}$, $\{(a_{n+1},\alpha_{n+1}),\ldots,(a_p,\alpha_p)\}$, 
$\{(b_1,\beta_1),\ldots,(b_m,\beta_m)\}$ and $\{(b_{m+1},\beta_{m+1}),\ldots,(b_q,\beta_q)\}$;  
\item[(ii)] If $n \geq 1$ and $q > m$, and if one of the $(a_i,\alpha_i)$ ($i=1,\ldots,n$) is equal
to one of the $(b_i,\beta_i)$ ($i=m+1,\ldots,q$)  then   
\begin{equation}
J^{m,n}_{p,q}\left[ z;\tau,\varepsilon \, \bigg| \, \begin{matrix} (a_i,\alpha_i)_{1,p} \\ (b_i,\beta_i)_{1,q-1}, 
(a_1,\alpha_1) \end{matrix}\right] = 
J^{m,n-1}_{p-1,q-1}\left[ z;\tau,\varepsilon \, \bigg| \, \begin{matrix} (a_i,\alpha_i)_{2,p} \\ (b_i,\beta_i)_{1,q-1} \end{matrix}\right] ;
\end{equation}
\item[(iii)] If $m \geq 1$ and $p > n$, and if one of the $(b_i,\beta_i)$ ($i=1,\ldots,m$) is equal
to one of the $(a_i,\alpha_i)$ ($i=n+1,\ldots,p$), then   
\begin{equation}
J^{m,n}_{p,q}\left[ z;\tau,\varepsilon \, \bigg| \, \begin{matrix} (a_i,\alpha_i)_{1,p-1},(b_1,\beta_1) \\ (b_i,\beta_i)_{1,q} \end{matrix}\right] = 
J^{m-1,n}_{p-1,q-1}\left[ z;\tau,\varepsilon \, \bigg| \, \begin{matrix} (a_i,\alpha_i)_{1,p-1} \\ (b_i,\beta_i)_{2,q} \end{matrix}\right] .
\end{equation}
\end{description}

\begin{prop}
\label{prop.2}
 The Fox-Barnes $J$-function satisfies the 
following properties: 
\begin{align}
& \textnormal{(i)} \; J^{m,n}_{p,q}\left[ z^{-1};\tau,\varepsilon \, \bigg| \, \begin{matrix} (a_i,\alpha_i)_{1,p}\\ (b_i,\beta_i)_{1,q} \end{matrix}\right]    = 
J^{n,m}_{q,p}\left[ z;\tau,\varepsilon \, \bigg| \, \begin{matrix} (1+\tau-b_i,\beta_i)_{1,q} \\ (1+\tau-a_i,\alpha_i)_{1,p}  \end{matrix}\right] , \\[1ex]
& \textnormal{(ii)} \; z^\sigma \, J^{m,n}_{p,q}\left[ z;\tau,\varepsilon \, \bigg| \, \begin{matrix} (a_i,\alpha_i)_{1,p} \\ (b_i,\beta_i)_{1,q}\end{matrix}\right] = 
\operatorname{e}^{\pi \varepsilon\sigma^2/\tau} \, J^{m,n}_{p,q}\left[\operatorname{e}^{-2\pi \varepsilon\sigma/\tau}  z;\tau,\varepsilon \, \bigg| \, \begin{matrix} (a_i+\sigma\alpha_i,\alpha_i)_{1,p}  \\ (b_i+\sigma\beta_i,\beta_i)_{1,q}  \end{matrix}\right] , \\[1ex]
& \textnormal{(iii)} \; \frac{1}{k} J^{m,n}_{p,q}\left[ z;\tau,\varepsilon \, \bigg| \, \begin{matrix} (a_i,\alpha_i)_{1,p} \\ (b_i,\beta_i)_{1,q}\end{matrix}\right]   = 
J^{m,n}_{p,q}\left[ z^k;\tau, k^2\varepsilon \, \bigg| \, \begin{matrix} (a_i,k\alpha_i)_{1,p} 
 \\ (b_i,k\beta_i)_{1,q} \end{matrix}\right]  , \quad (k > 0) . 
\end{align}
\end{prop}
\bigskip

\begin{proof} These properties follows directly from eq.\eqref{fox-barnes}
with a change of the integration variable, namely $s\to -s$ in (i), $s\to s+\sigma$ in (ii) and $s\to ks$ in (iii).  
\end{proof}

\bigskip

\begin{prop}
\label{prop.3} The Fox-Barnes $J$-function satisfies the 
following property:
\begin{equation}
 J^{m,n}_{p,q}\left[ z;\tau,\varepsilon \, \bigg| \, \begin{matrix} (a_i,\alpha_i)_{1,p} \\ (b_i,\beta_i)_{1,q} \end{matrix}\right]   = 
A \tau \, J^{m,n}_{p,q}\left[ (z/B)^\tau ;\tau^{-1},\varepsilon-\frac{\Delta_2 \, \log\tau}{2\pi} \, \bigg| \, \begin{matrix} 
(\tau^{-1}a_i,\alpha_i)_{1,p}   \\ (\tau^{-1}b_i,\beta_i)_{1,q}   \end{matrix}\right] ,
\end{equation} 
where 
\begin{equation}
\begin{aligned}
& A = (2\pi)^{\frac{(\tau-1)}{2\tau}[\delta^\ast + (1+\tau)(m+n-q)]} 
\tau^{\frac{1}{2\tau}[(1+\tau)D_1 - D_2] + \mathcal{N}} , \\[1ex]
& B = (2\pi)^{\frac{(\tau-1)}{2\tau}\Delta^\ast} \tau^{\frac{1}{2\tau}[(1+\tau)\Delta_1 - 2 \Pi]} , 
\end{aligned}
\end{equation}
with 
\begin{equation}
\begin{aligned}
& \delta^\ast = \sum_{i=1}^q b_i - \sum_{i=1}^p a_i , \qquad 
  \Delta^\ast = \sum_{i=1}^q \beta_i - \sum_{i=1}^p \alpha_i 
\end{aligned}
\end{equation}
and $\mathcal{N}$, $\Delta_i$, $D_i$, ($i=1,2$) and $\Pi$ defined by \eqref{def.N}, \eqref{def.Delta_k}, 
\eqref{def.D_k} and \eqref{def.Pi_1}, respectively. 
\end{prop}

\bigskip

\begin{proof} This property follows directly from eq.\eqref{fox-barnes} using
the modular property \eqref{modular} and a change of the integration variable $(s\to \tau s$). 
\end{proof}

\bigskip
\paragraph{Remark 2.}  Note that one of the effects of the modular transformation of the double gamma function in the definition of the Fox-Barnes function is to change the quantity $\varepsilon \to \varepsilon^\prime = \varepsilon-\Delta_2 \log{\tau} /(2\pi)$. Therefore, if in the definition of the Fox-Barnes function we had not used the term $\operatorname{e}^{\pi \varepsilon s^2/\tau}$, then due to this modular transformation we would have its equivalence with a function with the term $\operatorname{e}^{\pi \varepsilon^\prime s^2/\tau}$ with $\varepsilon^\prime \neq 0$. For this reason, in \cite{Kusnetsov} the term $\operatorname{e}^{\pi \varepsilon s^2/\tau}$ was introduced in the definition of the Meijer-Barnes function, and which we evidently also use in the definition of the Fox-Barnes function.  Note that the term $\operatorname{e}^{\pi \varepsilon s^2/\tau}$ 
was not needed in Examples 1, 3 and 4 because there we had $\tau = 1$, and so $\varepsilon^\prime = \varepsilon$ even if $\Delta_2 \neq 0$. Similarly for a hyper-balanced Fox-Barnes $J$-function this term is not needed too.

\bigskip

\begin{prop}
\label{prop.4} 
The Fox-Barnes $J$-function satisfies the 
following properties ($k~=~1,2,\ldots$):
\begingroup
\allowdisplaybreaks 
\begin{align}
&  \begin{aligned} \textnormal{(i)} \; & \frac{d^k\;}{dz^k}\; \left[ z^\omega J^{m,n}_{p,q}\left[ z^\sigma ;\tau,\varepsilon \, \bigg| \, \begin{matrix} (a_i,\alpha_i)_{1,p}\\ (b_i,\beta_i)_{1,q} \end{matrix}\right] \right] 
\\[1ex]
&  = 
\tau^k z^{\omega-k} J^{m,n+2}_{p+2,q+2}\left[ z^\sigma ;\tau,\varepsilon \, \bigg| \, \begin{matrix} 
(-\omega,\sigma),(-\omega+k+\tau,\sigma), (a_i,\alpha_i)_{1,p}\\ (b_i,\beta_i)_{1,q}, 
(-\omega+\tau,\sigma),(-\omega+k,\sigma) \end{matrix}\right] \\[1ex]
& = (-\tau)^k z^{\omega-k} J^{m+2,n}_{p+2,q+2}\left[ z^\sigma ;\tau,\varepsilon \, \bigg| \, \begin{matrix} (a_i,\alpha_i)_{1,p}, (1-\omega,\sigma), (k-1+\tau-\omega,\sigma) \\ 
(1+\tau-\omega,\sigma), (k-1-\omega,\sigma), (b_i,\beta_i)_{1,q} \end{matrix}\right] \\[1ex]
&  =  
\tau^k \operatorname{e}^{\pi\varepsilon[(\omega-k)/\sigma]^2/\tau}
J^{m,n+2}_{p+2,q+2}\left[ \operatorname{e}^{-2\pi\varepsilon(\omega-k)/\sigma}\,  z^\sigma ;\tau,\varepsilon \, \bigg| \, \begin{matrix} 
(-k,\sigma),(\tau,\sigma), (a_i+ \alpha_i(\omega-k)/\sigma,\alpha_i)_{1,p}\\ (b_i+\beta_i(\omega-k)/\sigma,\beta_i)_{1,q}, 
(\tau-k,\sigma),(0,\sigma) \end{matrix}\right] 
\end{aligned}\\[1ex] 
&  \; \begin{aligned}  \textnormal{(ii)} \; & \prod_{j=1}^k \left(z\frac{d\;}{dz}-c_j\right) \; \left[ z^\omega J^{m,n}_{p,q}\left[ z^\sigma ;\tau,\varepsilon \, \bigg| \, \begin{matrix} (a_i,\alpha_i)_{1,p}\\ (b_i,\beta_i)_{1,q} \end{matrix}\right] \right] 
\\[1ex]
&  = 
\tau^k z^{\omega} J^{m,n+2k}_{p+2k,q+2k}\left[ z^\sigma ;\tau,\varepsilon \, \bigg| \, \begin{matrix} 
(c_j-\omega,\sigma)_{1,k}, (1+\tau+c_j-\omega,\sigma)_{1,k}, (a_i,\alpha_i)_{1,p}\\ (b_i,\beta_i)_{1,q}, 
(\tau+c_j-\omega,\sigma)_{1,k}, (1+c_j-\omega,\sigma)_{1,k} \end{matrix}\right] \\[1ex]
&  = (-\tau)^k z^{\omega} J^{m+2k,n}_{p+2k,q+2k}\left[ z^\sigma ;\tau,\varepsilon \, \bigg| \, \begin{matrix} (a_i,\alpha_i)_{1,p}, (1+c_j-\omega,\sigma)_{1,k},(\tau+c_j-\omega,\sigma)_{1,k} \\ 
(1+\tau+c_j-\omega,\sigma)_{1,k}, (c_j-\omega,\sigma)_{1,k}, (b_i,\beta_i)_{1,q} \end{matrix}\right]  
\end{aligned}\\[1ex]
&  \; \begin{aligned}  \textnormal{(iii)} \; & \frac{d^k\;}{dz^k}\; \left[ J^{m,n}_{p,q}\left[ (cz+d)^\sigma ;\tau,\varepsilon \, \bigg| \, \begin{matrix} (a_i,\alpha_i)_{1,p}\\ (b_i,\beta_i)_{1,q} \end{matrix}\right] \right] 
\\[1ex]
&  = 
\frac{\tau^k c^k}{(cz+d)^k}   J^{m,n+2}_{p+2,q+2}\left[ (cz+d)^\sigma ;\tau,\varepsilon \, \bigg| \, \begin{matrix} 
(0,\sigma),(k+\tau,\sigma), (a_i,\alpha_i)_{1,p}\\ (b_i,\beta_i)_{1,q}, 
(\tau,\sigma),(k,\sigma) \end{matrix}\right]    \qquad (\sigma > 0) 
\end{aligned} \\[1ex]
&  \; \begin{aligned}  \textnormal{(iv)} \; & \frac{d^k\;}{dz^k}\; \left[ J^{m,n}_{p,q}\left[ (cz+d)^{-\sigma} ;\tau,\varepsilon \, \bigg| \, \begin{matrix} (a_i,\alpha_i)_{1,p}\\ (b_i,\beta_i)_{1,q} \end{matrix}\right] \right] 
\\[1ex]
&   = 
\frac{\tau^k c^k}{(cz+d)^k}   J^{m+2,n}_{p+2,q+2}\left[ (cz+d)^{-\sigma} ;\tau,\varepsilon \, \bigg| \, \begin{matrix} 
(1,\sigma),(1-k+\tau,\sigma), (a_i,\alpha_i)_{1,p}\\ (b_i,\beta_i)_{1,q}, 
(1+\tau,\sigma),(1-k,\sigma) \end{matrix}\right]  \qquad (\sigma > 0) 
\end{aligned} 
\end{align}
\endgroup
\end{prop}

\bigskip

\begin{proof} For example, in order to prove (i) we change the order of the derivative and the 
integration in the definition of the Fox-Barnes $J$-function, evaluate the derivative 
of $z^{\omega-\sigma s}$ and use eq.\eqref{ap.B.general.G.tau.2}. We have 
\begin{equation}
\begin{aligned}
\frac{d^k\;}{dz^k}\left(z^{\omega-\sigma s}\right) & = 
\frac{\Gamma(\omega-\sigma s + 1)}{\Gamma(\omega-\sigma s+1-k)} z^{\omega-k - \sigma s} \\[1ex]
& = 
\tau^k z^{\omega-k}z^{-\sigma s} 
\frac{G(1+\tau+\omega-\sigma s;\tau)G(1+\omega-k - \sigma s;\tau)}{G(1+\omega-\sigma s;\tau)G(1+\tau+\omega 
-k -\sigma s;\tau)} 
\end{aligned}
\end{equation}
which gives the first expression, or 
\begin{equation}
\begin{aligned}
\frac{d^k\;}{dz^k}\left(z^{\omega-\sigma s}\right) & = 
(-1)^k \frac{\Gamma(\omega+\sigma s+1)}{\Gamma(-\omega+\sigma s-1+k)} z^{\omega-k - \sigma s} \\[1ex]
& = 
(-1)^k\tau^k  z^{\omega-k}z^{-\sigma s} 
\frac{G(1+\tau-\omega+\sigma s;\tau)G(k-1-\omega + \sigma s;\tau)}{G(1-\omega+\sigma s;\tau)G(k-1+\tau-\omega 
+\sigma s;\tau)} 
\end{aligned}
\end{equation}
which gives the second expression. The third expression follows using the property (ii) in 
Proposition~\ref{prop.2} 
in the first expression. Properties (ii), (iii) and (iv) are proved analogously. 
\end{proof}

\bigskip

\begin{prop}
\label{prop.4.1}
The Fox-Barnes $J$-function satisfies the 
following properties:
\begingroup
\allowdisplaybreaks 
\begin{align}
&  \begin{aligned} \textnormal{(i)} \; \;
& z\frac{d\;}{dz} J^{m,n}_{p,q}\left[ z^\sigma ;\tau,\varepsilon \, \bigg| \, \begin{matrix} (a_i,\alpha_i)_{1,p}\\ (b_i,\beta_i)_{1,q} \end{matrix}\right]   = \frac{\sigma}{\alpha_1}(a_1-1) J^{m,n}_{p,q}\left[ z^\sigma ;\tau,\varepsilon \, \bigg| \, \begin{matrix} (a_i,\alpha_i)_{1,p}\\ (b_i,\beta_i)_{1,q} \end{matrix}\right] \\[1ex]
& \phantom{z\frac{d\;}{dz} J^{m,n}_{p,q} z^\sigma ;\tau,\varepsilon \, \bigg|} + \frac{\sigma \tau}{\alpha_1} J^{m,n+1}_{p+1,q+1}\left[ z^\sigma ;\tau,\varepsilon \, \bigg| \, \begin{matrix} (a_1+\tau,\alpha_1),(a_1-1,\alpha_1), (a_i,\alpha_i)_{2,p}\\ (b_i,\beta_i)_{1,q}, 
(a_1+\tau-1,\alpha_1) \end{matrix}\right] ,
\end{aligned} \\[1ex]
& \begin{aligned} \textnormal{(i$^\prime$)} \; \;
& z\frac{d\;}{dz} J^{m,n}_{p,q}\left[ z^\sigma ;\tau,\varepsilon \, \bigg| \, \begin{matrix} (a_i,\alpha_i)_{1,p}\\ (b_i,\beta_i)_{1,q} \end{matrix}\right]   = \frac{\sigma}{\alpha_1}(a_1-\tau) J^{m,n}_{p,q}\left[ z^\sigma ;\tau,\varepsilon \, \bigg| \, \begin{matrix} (a_i,\alpha_i)_{1,p}\\ (b_i,\beta_i)_{1,q} \end{matrix}\right] , \\[1ex]
& \phantom{z\frac{d\;}{dz} J^{m,n}_{p,q} z^\sigma ;\tau,\varepsilon \, \bigg|} + \frac{\sigma \tau}{\alpha_1} J^{m,n+1}_{p+1,q+1}\left[ z^\sigma ;\tau,\varepsilon \, \bigg| \, \begin{matrix} (a_1+1,\alpha_1),(a_1-\tau,\alpha_1), (a_i,\alpha_i)_{2,p}\\ (b_i,\beta_i)_{1,q}, 
(a_1+1-\tau,\alpha_1) \end{matrix}\right] ,
\end{aligned} \\[1ex]
&  \begin{aligned} \textnormal{(ii)} \; \;
& z\frac{d\;}{dz} J^{m,n}_{p,q}\left[ z^\sigma ;\tau,\varepsilon \, \bigg| \, \begin{matrix} (a_i,\alpha_i)_{1,p}\\ (b_i,\beta_i)_{1,q} \end{matrix}\right]   = \frac{\sigma}{\beta_1}(b_1-1) J^{m,n}_{p,q}\left[ z^\sigma ;\tau,\varepsilon \, \bigg| \, \begin{matrix} (a_i,\alpha_i)_{1,p}\\ (b_i,\beta_i)_{1,q} \end{matrix}\right] \\[1ex]
& \phantom{z\frac{d\;}{dz} J^{m,n}_{p,q} z^\sigma ;\tau,\varepsilon \, \bigg|} + \frac{\sigma \tau}{\beta_1} J^{m+1,n}_{p+1,q+1}\left[ z^\sigma ;\tau,\varepsilon \, \bigg| \, \begin{matrix}  (a_i,\alpha_i)_{1,p}, (b_1+\tau-1,\beta_1)\\ (b_1-1,\beta_1), (b_1+\tau,\beta_1), (b_i,\beta_i)_{2,q} \end{matrix}\right] ,
\end{aligned} \\[1ex] 
& \begin{aligned} \textnormal{(ii$^\prime$)} \; \;
& z\frac{d\;}{dz} J^{m,n}_{p,q}\left[ z^\sigma ;\tau,\varepsilon \, \bigg| \, \begin{matrix} (a_i,\alpha_i)_{1,p}\\ (b_i,\beta_i)_{1,q} \end{matrix}\right]   = \frac{\sigma}{\beta_1}(b_1-\tau) J^{m,n}_{p,q}\left[ z^\sigma ;\tau,\varepsilon \, \bigg| \, \begin{matrix} (a_i,\alpha_i)_{1,p}\\ (b_i,\beta_i)_{1,q} \end{matrix}\right] \\[1ex]
& \phantom{z\frac{d\;}{dz} J^{m,n}_{p,q} z^\sigma ;\tau,\varepsilon \, \bigg|} + \frac{\sigma \tau}{\beta_1} J^{m+1,n}_{p+1,q+1}\left[ z^\sigma ;\tau,\varepsilon \, \bigg| \, \begin{matrix}  (a_i,\alpha_i)_{1,p}, (b_1+1-\tau,\beta_1)\\ (b_1-\tau,\beta_1), (b_1+1,\beta_1), (b_i,\beta_i)_{2,q} \end{matrix}\right] ,
\end{aligned} \\[1ex]
&  \begin{aligned} \textnormal{(iii)} \; \;
& z\frac{d\;}{dz} J^{m,n}_{p,q}\left[ z^\sigma ;\tau,\varepsilon \, \bigg| \, \begin{matrix} (a_i,\alpha_i)_{1,p}\\ (b_i,\beta_i)_{1,q} \end{matrix}\right]   = \frac{\sigma}{\beta_q}b_q J^{m,n}_{p,q}\left[ z^\sigma ;\tau,\varepsilon \, \bigg| \, \begin{matrix} (a_i,\alpha_i)_{1,p}\\ (b_i,\beta_i)_{1,q} \end{matrix}\right] \\[1ex]
& \phantom{z\frac{d\;}{dz} J^{m,n}_{p,q} z^\sigma ;\tau,\varepsilon \, \bigg|} + \frac{\sigma \tau}{\beta_q} J^{m,n+1}_{p+1,q+1}\left[ z^\sigma ;\tau,\varepsilon \, \bigg| \, \begin{matrix} (1+\tau+b_q,\beta_1), (a_i,\alpha_i)_{1,p}\\ (b_i,\beta_i)_{1,q-1}, 
(b_q+\tau,\beta_q), (b_q+1,\beta_q) \end{matrix}\right] ,
\end{aligned} \\[1ex]
&  \begin{aligned} \textnormal{(iv)} \; \;
& z\frac{d\;}{dz} J^{m,n}_{p,q}\left[ z^\sigma ;\tau,\varepsilon \, \bigg| \, \begin{matrix} (a_i,\alpha_i)_{1,p}\\ (b_i,\beta_i)_{1,q} \end{matrix}\right]   = \frac{\sigma}{\alpha_p}(a+p-1-\tau) J^{m,n}_{p,q}\left[ z^\sigma ;\tau,\varepsilon \, \bigg| \, \begin{matrix} (a_i,\alpha_i)_{1,p}\\ (b_i,\beta_i)_{1,q} \end{matrix}\right] \\[1ex]
& \phantom{z\frac{d\;}{dz} J^{m,n}_{p,q} z^\sigma ;\tau,\varepsilon \, \bigg|} - \frac{\sigma \tau}{\alpha_p} J^{m+1,n}_{p+1,q+1}\left[ z^\sigma ;\tau,\varepsilon \, \bigg| \, \begin{matrix} (a_1,\alpha_i)_{1,p-1}, (a_p-\tau,\alpha_p), (a_p-1,\alpha_p)\\ (a_p-1-\tau,\alpha_p), (b_i,\beta_i)_{1,q} \end{matrix}\right] ,
\end{aligned} 
\end{align}
\endgroup
where $n \geq 1$ in \textnormal{(i)} and \textnormal{(i$^\prime$)}, $m \geq 1$ in \textnormal{(ii)} and \textnormal{(ii$^\prime$)}, $q \geq m + 1$ in \textnormal{(iii)} and $p \geq n + 1$ in~\textnormal{(iv)}. 
\end{prop}

\bigskip

\begin{proof} Firstly we prove some identities. 
Taking $z\to z-1$ in eq.\eqref{ap.B.general.G.tau}, multiplying 
the result by $(z-1)/\tau$, using the functional relation of the gamma function and again 
eq.\eqref{ap.B.general.G.tau} we obtain the relation 
\begin{equation}
\label{prop.4.1.aux.1}
(z-1)G(z,\tau) = \tau \frac{G(z-1;\tau) G(z+\tau;\tau)}{G(z-1+\tau;\tau)} . 
\end{equation}
Taking $z\to z-1$ in eq.\eqref{ap.B.general.G.tau.2}, multiplying 
the result by $z-\tau$, using the functional relation of the gamma function and again 
eq.\eqref{ap.B.general.G.tau.2} we obtain the relation 
\begin{equation}
\label{prop.4.1.aux.2}
(z-\tau)G(z,\tau) = \tau \frac{G(z+1;\tau) G(z-\tau;\tau)}{G(z+1-\tau;\tau)} . 
\end{equation}
Writing eq.\eqref{ap.B.general.G.tau} as 
\begin{equation}
\frac{1}{G(z;\tau)} = \frac{\Gamma(z/\tau)}{\Gamma(z+1;\tau)} , 
\end{equation}
multiplying it by $z/\tau$, using the functional relation of the gamma function and
again eq.\eqref{ap.B.general.G.tau} we obtain  the relation 
\begin{equation}
\label{prop.4.1.aux.3}
\frac{z}{G(z;\tau)} = \tau \frac{G(z+1+\tau;\tau)}{G(z+1;\tau) G(z+\tau;\tau)} . 
\end{equation}
However, unlike in eq.\eqref{prop.4.1.aux.1} and eq.\eqref{prop.4.1.aux.2}, 
if we use eq.\eqref{ap.B.general.G.tau.2} to write an expression for $1/G(z;\tau)$ and 
repeat the above steps we do not obtain a new relation but we obtain again eq.\eqref{prop.4.1.aux.3}.
Then, in order to prove (i) and (ii), we evaluate the action of $z\frac{d\;}{dz}$ on the definition of the
Fox-Barnes $J$-function and use the relation in eq.\eqref{prop.4.1.aux.1}. In order to 
prove (i$^\prime$) and (ii$^\prime$), we again evaluate the action of $z\frac{d\;}{dz}$ on the definition of the
Fox-Barnes $J$-function and in this case we use the relation in eq.\eqref{prop.4.1.aux.2}. 
Finally, in order to prove (iii) and (iv), we repeat the same steps and use 
eq.\eqref{prop.4.1.aux.3}. 
\end{proof}

\section{Application: Laplace Transform of the KS function}
\label{sec.6}

The motivation for evaluating the Laplace transform of the Kilbas-Saigo (KS) function 
arises naturally from the need to analyse and characterize stochastic models 
governed by stretched fractional dynamics. Indeed in \cite{BLV} the authors 
extend classical renewal processes by introducing a non-local operator 
$\mathcal{D}^{(\alpha,\gamma)}_t$, whose solution involves the KS function. 
This function defines the survival probabilities of interarrival times and 
it plays a key role in how the related counting process behaves. To compute some 
essential statistical quantities like variances and correlations of the
process, it becomes essential to evaluate the Laplace transform of the KS function. 
Moreover, it is also crucial for expressing the model in terms of a time-changed representation, 
where generalized renewal processes are linked to subordinated Poisson processes.  

The KS function, denoted by $E_{a,m,l}(z)$, is defined as \cite{Simon,CMV}
\begin{equation}
\label{def.KS.function}
E_{a,m,l}(z) = \sum_{n=0}^\infty c_n z^n ,  
\end{equation}
with 
\begin{equation}
\label{def.coeff.KS.function}
c_0 = 1 , \qquad c_n = \prod_{k=0}^{n-1} \frac{\Gamma[1+a(km+l)]}{\Gamma[1+a(km+l+1)]} , 
\end{equation}
and where the parameters are such that $\alpha, m > 0$ and $l > -1/a$. 
This series converges for $z \in \mathbb{C}$. 
The fractional differential equation 
\begin{equation}
\label{s.diff.op}
\mathcal{D}_t^{(\alpha,\gamma)} f(t) + \lambda f(t) = 0 ,
\end{equation}
with 
\begin{equation}
\label{eq.1}
\operatorname{\mathcal{D}}^{(\alpha,\gamma)}_t = t^{-\gamma} {\sideset{_{\scriptscriptstyle C}^{}}{_t^{(\alpha)}}{\operatorname{\mathcal{D}}}} ,
\end{equation}
where ${\sideset{_{\scriptscriptstyle C}^{}}{_t^{(\alpha)}}{\operatorname{\mathcal{D}}}} $ is the
Caputo derivative of order $\alpha$ (for $0 < \alpha < 1$) with respect to $t$ \cite{Podlubny,Kilbas}, has solution 
\begin{equation}\label{KS}
f(t) = f(0) \, \operatorname{E}_{a,m,l}(-\lambda t^\nu) 
\end{equation}
with 
\begin{equation}
a = \alpha, \qquad m = 1+\gamma/\alpha, \qquad l = \gamma/\alpha, \qquad \nu = \alpha + \gamma. 
\end{equation}

Let $F(z) = \mathcal{L}_t[f(t);z]$ be the Laplace transform of $f(t)$ with respect to $t$ and $z$ is the Laplace
dual variable. Assuming we can evaluate the Laplace transform of $\operatorname{E}_{a,m,l}(-\lambda t^\nu)$ using the term-by-term evaluation of eq.\eqref{def.KS.function}, we obtain 
\begin{equation}
\label{laplace.1}
\mathcal{L}_t[\operatorname{E}_{a,m,l}(-\lambda t^\nu);z] = \sum_{n=0}^\infty 
c_n (-\lambda)^n \frac{\Gamma(1+n \nu)}{z^{n\nu+1}} . 
\end{equation}
Using the ratio test and Stirling approximation we see that if $\nu < a$ (i.e., $\gamma < 0$) 
this series \textit{converges} and if 
$\nu > a$ (i.e., $\gamma > 0$) it \textit{diverges}. 

There is interest in the case $\gamma > 0$. Indeed, in \cite{BLV} it is show 
that when $\alpha + \gamma > 1$ (which requires $\gamma > 1-\alpha > 0$) the 
expected value of the interarrival time of the renewal process 
defined from eq.\eqref{eq.1} is finite. Thus we are interested in evaluating
the Laplace transform of the KS function when $\nu > a$. This question 
was addressed in \cite{BLV} and in this section we want to show how it 
can be expressed in terms of the Fox-Barnes $J$-function. 

In what follows we will consider Fox-Barnes functions defined using the contour $\mathcal{C}_{ic\infty}$. Thus the parameters must be such that the cases (a), (c), (o) and (q) discussed in Subsection~\ref{section.3.1} hold. For the contour $\mathcal{C}_{ic\infty}$ we also have $\mathcal{A} < c < \mathcal{B}$. 

We want to show the following 
\begin{prop}
\label{prop.5}
Let $J^{m,n}_{p,q}\left[ \lambda t^\nu ;\tau,\varepsilon \, \big| \, \begin{matrix} (a_i,\alpha_i)_{1,p}; (b_i,\beta_i)_{1,q} \end{matrix}\right]$ be a Fox-Barnes $J$-function satisfying a VL condition. Let $\mathcal{A}$ and $\mathcal{B}$ as in eq.\eqref{def.A.B} and 
assume $\mathcal{A} < \mathcal{B}$, $\mathcal{A} < 1/\nu$ ($\nu>0$) and 
in addition $\mathcal{A} < c < \mathcal{B}$ in the case of the VL condition~4. Then its Laplace transform  exists 
and is given by 
\begin{equation}
\label{Laplace.I}
\begin{aligned}
\mathcal{L}_t\left[ J^{m,n}_{p,q}\left[ \lambda t^\nu ;\tau,\varepsilon \, \bigg| \, \begin{matrix} (a_i,\alpha_i)_{1,p}\\ (b_i,\beta_i)_{1,q} \end{matrix}\right]; z\right]  = 
\frac{\tau^{1/2}}{(2\pi)^{(\tau-1)/2}} \, \frac{1}{z}  
J^{m,n+1}_{p+1,q+1}\left[ \frac{\lambda \tau^\nu}{z^\nu} ;\tau,\varepsilon \, \bigg| \, \begin{matrix} (0,\nu), (a_i,\alpha_i)_{1,p}\\ (b_i,\beta_i)_{1,q}, (\tau,\nu) \end{matrix}\right] & \\[1ex]
 =\frac{\tau^{1/2}}{(2\pi)^{(\tau-1)/2}} \, \frac{1}{z}  
J^{n+1,m}_{q+1,p+1}\left[ \frac{z^\nu}{\lambda \tau^\nu} ;\tau,\varepsilon \, \bigg| \, \begin{matrix} (1+\tau-b_i,\beta_i)_{1,q}, (1,\nu)\\(1+\tau,\nu), (1+\tau-a_i,\alpha_i)_{1,p} \end{matrix}\right] &
\end{aligned}
\end{equation}
for $z \in \mathbb{C}$ and $\operatorname{Re}z > 0$. 
\end{prop}

\bigskip

\begin{proof} Under the above assumptions, the function $J^{m,n}_{p,q}\left[ \lambda t^\nu ;\tau,\varepsilon \, \big| \, \begin{matrix} (a_i,\alpha_i)_{1,p}; (b_i,\beta_i)_{1,q} \end{matrix}\right]$ defined in terms of the
contour $\mathcal{C}_{ic\infty}$ exists. Indeed in the case of 
the VL conditions 1, 2 and 3 (see Definition~\ref{def.1}) we 
choose $c$ such that $\mathcal{A} < c < \mathcal{B}$, 
while in case 4 we have $\mathcal{A} < c < \mathcal{B}$ by assumption. In order to evaluate the Laplace transform, we will assume we can change the order of integration, which gives the integral $\int_0^\infty \operatorname{e}^{-zt} t^{-\nu s} dt = \Gamma(1-\nu s) z^{\nu s}$ if $1-\nu \,\operatorname{Re} s > 0$. 
This is justified                  
for $f(t) = J_{p,q}^{m,n}[\lambda t^\mu; \cdot]$ because of the absolute 
convergence of the double integral (Fubini-Tonelli theorem) which follows 
from the absolute convergence of the Mellin transform of $J_{p,q}^{m,n}[\lambda t^\mu; \cdot]$ with a VL condition together with the 
exponential damping of the Laplace kernel.
Thus if $1/\nu > \mathcal{A}$, we can choose $c$ such that $\mathcal{A} < c < \operatorname{min}(\mathcal{B},1/\nu)$. Next we use eq.\eqref{ap.B.general.G.tau.2} to write $\Gamma(1-\nu s)$ in terms of $G(1+\tau -\nu s;\tau)$ and $G(1-\nu s;\tau)$ and the definition \eqref{fox-barnes} to identify the resulting expression as the RHS in eq.\eqref{Laplace.I}. The two expressions are related by property 
(i) in Proposition~\ref{prop.2}. Now let $\Delta_2^{\scriptscriptstyle L}$, $\Theta_2^{\scriptscriptstyle L}$, etc, the quantities analogous to $\Delta_2$, $\Theta_2$, etc, but defined using the parameters of the function $J^{n+1,m}_{q+1,p+1}\left[ \frac{z^\nu}{\lambda \tau^\nu} ;\tau,\varepsilon \, \big| \, \begin{matrix} (1+\tau-b_i,\beta_i)_{1,q}, (1,\nu); (1+\tau,\nu), (1+\tau-a_i,\alpha_i)_{1,p} \end{matrix}\right]$. These quantities are expected to be related and indeed it is not difficult to see that 
\begin{equation}
\begin{split}
& \Delta_2^{\scriptscriptstyle L} = \Delta_2, \qquad \Theta_2^{\scriptscriptstyle L} = \Theta_2, \qquad 
\Omega_2^{\scriptscriptstyle L} = \Omega_2^\dagger, \qquad 
\Omega^\dagger{}^{\scriptscriptstyle L} = \Omega_2 , \qquad 
\Delta^\ast{}^{\scriptscriptstyle L} = - \Delta\\
& \Pi^{\scriptscriptstyle L} = -\Pi + (1+\tau)\Delta_1 + \tau \nu, \qquad \Upsilon^{\scriptscriptstyle L} = -(1+\tau)\Omega_1^\dagger -\Upsilon^\dagger \qquad  \Upsilon^\dagger{}^{\scriptscriptstyle L} = 
-\tau \nu - (1+\tau) \Omega_1 - \Upsilon \\
& \Lambda^\dagger = -\Lambda + \tau \nu \log\nu + (1+\tau)\Theta_1,  \qquad 
D_1^{\scriptscriptstyle L} = -D_1 + (1+\tau)\mathcal{N} +\tau, \\
&  D_2^{\scriptscriptstyle L} = D_2 - 2(1+\tau) D_1 + (1+\tau)^2 \mathcal{N} + \tau(\tau+2) , \\
& \Phi_1^{\scriptscriptstyle L} + 2 \Phi_3^{\scriptscriptstyle L} = 
\Phi_1 + 2\Phi_3 + \nu + \frac{2c}{\tau}(\Omega_2 -\Omega_2^\dagger) . 
\end{split}
\end{equation}
It is not difficult to see that if the parameters of the original function $J^{m,n}_{p,q}[\cdot]$ satisfy the VL conditions, then the parameters of its Laplace transform $J^{n+1,m}_{q+1,p+1}[\cdot]$ also satisfy the VL conditions. For conditions 1 and 2 it follows directly from the first two expressions above. For conditions 3 and 4, we notice that since $\Delta_2 = 0$ and $\Omega_2 = 0$, it follows from $\Delta_2 = \Omega_2 + \Omega_2^\dagger$ that also $\Omega_2^\dagger = 0$, and so from the last expression above we have $\Phi_1^{\scriptscriptstyle L} + 2 \Phi_3^{\scriptscriptstyle L} = \Phi_1 + 2\Phi_3 + \nu$. So we have $\Phi_1^\dagger + 2\Phi_3^\dagger >0$ since $\nu > 0$ and $\Phi_1 + 2 \Phi_3 \geq 0$, and then we can use the contour $\mathcal{C}_{ic\infty}$ for $|\operatorname{arg}z - \operatorname{Im}\Lambda| < (\Phi_1 + 2\Phi_3 + \nu)/2$ even in the case of VL condition 4 for $J^{m,n}_{p,q}[\cdot]$. 
Finally, the poles of $J^{n+1,m}_{q+1,p+1}[\cdot]$ are 
$\underline{\acute{\mys}}_{m,n}^{i} = (a_j + m\tau + n)/\alpha_j$ for  
$m,n = 0,1,2,\ldots$, $j = n+1,\ldots, p$, 
$\underline{\grave{\mys}}_{m^\prime,n^\prime}^i = 
- (-b_i + m^\prime \tau + n^\prime)/\beta_i$ for  
$m^\prime, n^\prime = 1,2,\ldots$, $i = m+1,\ldots, q$, and 
$\underline{\grave{\mys}}_{m,n}^\nu = 
-(1+m\tau+n)/\nu$, $m,n = 0,1,2,\ldots$. 
Let $\mathcal{A}^\prime$ and $\mathcal{B}^\prime$ be defined as 
\begin{equation}
\mathcal{A}^\prime = \underset{n+1 \leq j \leq p}{\operatorname{min}} \;
\frac{\operatorname{Re} a_j}{\alpha_j} , \qquad \mathcal{B}^\prime = - 
\underset{m+1\leq i \leq q}{\operatorname{min}}\frac{1}{\nu}(1 + \tau- \operatorname{Re} b_i) . 
\end{equation}
Then we must have for $\operatorname{Re}s = c^\prime$ along the vertical contour that 
\begin{equation}
\operatorname{max}(\mathcal{B}^\prime,-1/\nu) < c^\prime < 
\mathcal{A}^\prime . 
\end{equation}
Let us recall that property (i) in Proposition~\ref{prop.2} was proved using the change of variable $s\to -s$. Thus $c^\prime = -c$. 
Since $\mathcal{A}^\prime = -\mathcal{A}$ and $\mathcal{B}^\prime = -
\mathcal{B}$ we see that $c^\prime$ indeed satisfy the above 
inequalities because 
$\mathcal{A} < c < \operatorname{min}(\mathcal{B},1/\nu)$. 
We conclude therefore 
that the function $J^{n+1,m}_{q+1,p+1}[\cdot]$ defined in terms of the contour $\mathcal{C}_{ic^\prime\infty}$ exists. 
\end{proof}

\medskip

\enlargethispage{2ex}

Next we will show that the function 
$E_{a,m,l}\left(-\lambda t^\nu\right)$ can be written in the form of a Fox-Barnes $J$-function. 

\begin{prop}
\label{prop.6}
Let $\varphi = (1+ al)\tau$ and $\tau = 1/(am)$. Then 
\begin{equation}
\label{eq.prop.6}
\begin{aligned}
E_{a,m,l}\left(-\lambda t^\nu\right) &  = 
\frac{G(\varphi+a\tau;\tau)}{(2\pi)^{\tau-1}\nu G(\varphi;\tau)} 
 J^{1,2}_{3,3}\left[ \lambda^{1/\nu}t ;\tau  \, \bigg| \, \begin{matrix} (0,1/\nu), (1+\tau-\varphi,1/\nu), (0,1/\nu) \\ (\tau,1/\nu), (\tau,1/\nu), (1+\tau-\varphi-a\tau,1/\nu) \end{matrix}\right] \\[1ex]
 & = \frac{G(\varphi+a\tau;\tau)}{(2\pi)^{\tau-1} G(\varphi;\tau)} 
 J^{1,2}_{3,3}\left[ \lambda t^\nu ;\tau  \, \bigg| \, \begin{matrix} (0,1), (1+\tau-\varphi,1), (0,1) \\ (\tau,1), (\tau,1), (1+\tau-\varphi-a\tau,1) \end{matrix}\right] .
 \end{aligned}
\end{equation}
\end{prop}

\bigskip

\begin{proof} Firstly, let us write the explicit 
expression for this particular Fox-Barnes $J$-function, that is, 
\begin{equation}
\begin{split}
& \frac{G(\varphi+a\tau;\tau)}{(2\pi)^{\tau-1}\nu G(\varphi;\tau)} J^{1,2}_{3,3}\left[ \lambda^{1/\nu}t ;\tau  \, \bigg| \, \begin{matrix} (0,1/\nu), (1+\tau-\varphi,1/\nu), (0,1/\nu) \\ (\tau,1/\nu), (\tau,1/\nu), (1+\tau-\varphi-a\tau,1/\nu) \end{matrix}\right] \\[1ex]
 & \quad = \frac{G(\varphi+a\tau;\tau)}{(2\pi)^{\tau-1}\nu G(\varphi;\tau)} \frac{1}{2\pi i}\int_{\mathcal{C}} (\lambda^{1/\nu}t)^{-s} 
 \frac{G(\tau+s/\nu;\tau) G(\varphi-s/\nu;\tau) G(1+\tau-s/\nu;\tau)}{G(1-s/\nu;\tau) G(\varphi+a\tau-s/\nu;\tau) G(s/\nu;\tau)}\, ds .
 \end{split}
\end{equation}
In relation to its existence, note that we have 
\begin{equation}
\Delta_2 = \Delta_1 = 0, \quad \Theta_2 = \Theta_1 = 0, \quad \Omega_2 = \Omega_1 = 0, \quad \Delta^\ast = 0, \quad 
\epsilon = 0, 
\end{equation}
and so all cases discussed in (i) and (ii) in Subsection~\ref{section.3.1} are ruled out. Next we consider 
the cases in (iii) and we have 
\begin{equation}
\Pi = \frac{a\tau}{\nu} 
\end{equation}
and then 
\begin{equation}
\Phi_1 = \frac{a}{\nu}.
\end{equation}
Thus $\Phi_1 > 0$ and we can use the contour $\mathcal{C}_{-\infty}$, 
and since $\operatorname{Im}\Pi = 0$ the contour $\mathcal{C}_{-\infty}$ 
can be deformed into a contour $\mathcal{C}_{(\theta_-,\theta_+)}$ 
with $-\pi < \theta_- < -\pi/2$ and $\pi/2 < \theta_+ < \pi$. 
In relation to the contour $\mathcal{C}_{ic\infty}$, since $\operatorname{Im} \Pi = 0$ we have $K_3 = 0$ and we go to the 
cases (iv) related to $\mathcal{C}_{ic\infty}$. We also have 
\begin{equation}
\Upsilon = \frac{\tau}{\nu}(1-a), \qquad \Lambda = -\frac{a\tau}{\nu} \log\nu , 
\end{equation} 
which gives
\begin{equation}
\Phi_3 = \frac{1-a}{\nu} , \qquad \operatorname{Im}\Lambda = 0 , 
\end{equation}
and then
\begin{equation}
\Phi_1 + 2 \Phi_3 = \frac{1}{\nu} . 
\end{equation}
Thus we can also use the contour $\mathcal{C}_{ic\infty}$ for 
\begin{equation}
|\operatorname{arg}t| < \frac{\pi}{2\nu} . 
\end{equation}
The poles of $J^{1,2}_{3,3}[\cdot]$  are 
$\grave{\mys}_{m,n} = -\nu(m\tau+n)$ from 
$G(s/\nu;\tau)$, $\acute{\mys}^{(1)}_{m,n} = \nu(1+m\tau+n)$ from 
$G(1-s/\nu;\tau)$ and $\acute{\mys}^{(2)}_{m,n} = \nu(\varphi+a\tau+m\tau+n)$ from $G(\varphi+a\tau-s/\nu;\tau)$ (where $m,n = 0,1,2,\ldots$). So we have $\mathcal{A} = 0$, $\mathcal{B} = \operatorname{min}(\nu,\nu(\varphi+a\tau)) > 0$ and $\mathcal{A} < c < \mathcal{B}$.

The contour $\mathcal{C}_{ic\infty}$ can therefore be deformed into $\mathcal{C}_{-\infty}$ and the residue theorem can be used to obtain a series representation for the RHS of eq.\eqref{eq.prop.6}. The poles we have to consider are $\grave{\mys}_{m,n} = -\nu(m\tau+n)$ and the order of them depends on the value of $\tau$. However, for this specific $J$-function, there is a function $G(\tau+s/\nu;\tau)$ in the numerator and $G(\tau+\grave{\mys}_{m,n}/\nu;\tau) = 0$ for $m=1,2,\ldots$ and $n=0,1,\ldots$. This expected since $G(\tau+s/\nu;\tau)$ and $G(s/\nu;\tau)$ are related by eq.\eqref{ap.B.general.G.tau.2}. The same happens for 
$G(1+\tau-s/\nu;\tau)$ and $G(1-s/\nu;\tau)$. Thus using eq.\eqref{ap.B.general.G.tau.2} we can write 
\begin{equation}
E_{a,m,l}\left(-\lambda t^\nu\right) = \frac{G(\varphi+a\tau;\tau)}{\nu G(\varphi;\tau)} \frac{1}{2\pi i}\int_{\mathcal{C}} (\lambda^{1/\nu}t)^{-s} 
 \frac{\Gamma(s/\nu) \Gamma(1-s/\nu)G(\varphi-s/\nu;\tau)}{G(\varphi+a\tau-s/\nu;\tau)}\, ds 
\end{equation}
Using the residue theorem taking into account the 
poles of the gamma function $\Gamma(s/\nu)$ at $\grave{\mys}_n = - n \nu$ we obtain 
\begin{equation}
E_{a,m,l}\left(-\lambda t^\nu\right) = \frac{G(\varphi+a\tau;\tau)}{ G(\varphi;\tau)} \sum_{n=0}^\infty  (\lambda^{1/\nu}t)^{n\nu} 
\frac{(-1)^n}{n!} \frac{\Gamma(1+n)G(\varphi+n;\tau)}{G(\varphi+a\tau+n;\tau)}
\end{equation}
Finally using eq.\eqref{ap.B.general.G.tau.n} we obtain for $c_n$ 
in eq.\eqref{def.coeff.KS.function} that 
\begin{equation}
c_n = \frac{G(\varphi+a\tau;\tau) G(\varphi+n;\tau)}{G(\varphi+a\tau+n;\tau) G(\varphi;\tau)} ,
\end{equation}
and this gives the series representation of $E_{a,m,l}\left(-\lambda t^\nu\right)$.
\end{proof}

\medskip

Note that the second expression in eq.\eqref{eq.prop.6} is obtained 
from the first one using the property in Proposition~\ref{prop.2}-(iii). 
Since $\alpha_i = 1$ and $\beta_i = 1$ ($i=1,2,3$) this is an
example of a Meijer-Barnes $K$-function introduced in 
\cite{Kusnetsov}. Its Laplace transform however cannot be
written as this particular case, as we will see. 

\begin{prop}
\label{prop.7}
The Laplace transform of the KS function $E_{a,m,l}\left(-\lambda t^\nu\right)$ is given by 
\begin{equation}
\begin{aligned}
& \mathcal{L}_t[ E_{a,m,l}\left(-\lambda t^\nu\right);z] \\[1ex]
& \qquad = 
\frac{\tau^{1/2}G(\varphi+a\tau;\tau)}{(2\pi)^{3(\tau-1)/2}G(\varphi;\tau)} \,\frac{1}{z} \, 
 J^{1,3}_{4,4}\left[ \frac{\lambda \tau^\nu}{z^\nu} ;\tau  \, \bigg| \, \begin{matrix} (0,\nu), (0,1), (1+\tau-\varphi,1), (0,1) \\ (\tau,1), (1+\tau-\varphi-a\tau,1), (\tau,1), (\tau,\nu) \end{matrix}\right] \\[1ex]
 & \qquad = 
\frac{\tau^{1/2}G(\varphi+a\tau;\tau)}{(2\pi)^{3(\tau-1)/2}G(\varphi;\tau)} \,\frac{1}{z} \, 
 J^{3,1}_{4,4}\left[ \frac{ z^\nu}{\lambda \tau^\nu} ;\tau  \, \bigg| \, \begin{matrix} (1,1), (\varphi+a\tau,1), (1,1), (1,\nu)\\ (1+\tau,\nu), (1+\tau,1), (\varphi,1), (1+\tau,1) \end{matrix}\right] 
\end{aligned}
\end{equation}

\end{prop}

\bigskip

\begin{proof} This is a direct consequence of propositions \ref{prop.5} and \ref{prop.6} 
and the property in Proposition \ref{prop.2}-(i). The explicit expression we have is 
\begin{equation}
\begin{aligned}
& J^{3,1}_{4,4}\left[ \frac{ z^\nu}{\lambda \tau^\nu} ;\tau  \, \bigg| \, \begin{matrix} (1,1), (\varphi+a\tau,1), (1,1), (1,\nu)\\ (1+\tau,\nu), (1+\tau,1), (\varphi,1), (1+\tau,1) \end{matrix}\right] \\[1ex]
& \frac{1}{2\pi i}\int_{\mathcal{C}} \left(\frac{z^\nu}{\lambda \tau^\nu}\right)^{-s} \frac{G(1+\tau+\nu s;\tau) G(1+\tau+s;\tau) G(\varphi+s;\tau) G(\tau-s;\tau)}{G(-s;\tau)G(\varphi+a\tau+s;\tau) G(1+s;\tau) G(1+\nu s;\tau)} ds . 
\end{aligned}
\end{equation}
In relation to the existence of this particular 
Fox-Barnes $J$-function, note that we have 
\begin{equation}
\Delta_2^{\scriptscriptstyle L} = \Delta_1^{\scriptscriptstyle L} = 0, \quad \Theta_2^{\scriptscriptstyle L} = \Theta_1^{\scriptscriptstyle L} = 0, \quad \Omega_2^{\scriptscriptstyle L} = \Omega_1^{\scriptscriptstyle L} = 0, \quad \Delta^\ast{}^{\scriptscriptstyle L} = 0, \quad 
\epsilon = 0, 
\end{equation}
where we used the notation introduced in the proof of Proposition~\ref{prop.5}. 
Thus all cases discussed in (i) and (ii) in Subsection~\ref{section.3.1} are ruled out. Next we consider 
the cases in (iii) and we have 
\begin{equation}
\Pi^{\scriptscriptstyle L} = \tau(\nu-a)  
\end{equation}
and then 
\begin{equation}
\Phi_1^{\scriptscriptstyle L} = \nu - a .
\end{equation}
So if $\nu > a$ we can use the contour $\mathcal{C}_{-\infty}$ (case (i)) and if 
$\nu < a$ we can use the contour $\mathcal{C}_{+\infty}$ (case (j)). 
Moreover, since $\operatorname{Im}\Pi^\dagger = 0$ the contour $\mathcal{C}_{-\infty}$ 
can be deformed into a contour $\mathcal{C}_{(\theta_-,\theta_+)}$ 
with $-\pi < \theta_- < -\pi/2$ and $\pi/2 < \theta_+ < \pi$ and 
the contour $\mathcal{C}_{+\infty}$ 
can be deformed into a contour $\mathcal{C}_{(\theta_-,\theta_+)}$ 
with $-\pi/2 < \theta_- < 0$ and $0 < \theta_+ < \pi/2$. 
As discussed in (iv), we can also use the contour $\mathcal{C}_{ic\infty}$ (case (o)) 
since in this case $K_3 = 0$. We have 
\begin{equation}
\Upsilon^{\scriptscriptstyle L} = \tau , \quad \Lambda^{\scriptscriptstyle L}= \tau \nu \log\nu .
\end{equation}
Then $\Phi_3^{\scriptscriptstyle L} = 1$ and 
\begin{equation}
\Phi_1^{\scriptscriptstyle L} + 2 \Phi_3^{\scriptscriptstyle L} 
= \nu - a + 2 . 
\end{equation}
We also recall that here we have $z^\nu$ instead of $z$. Thus in the 
condition for the case (o) we must use 
$\nu |\operatorname{arg}z|$ instead of $|\operatorname{arg}z|$. 
Consequently we can use $\mathcal{C}_{ic\infty}$ if 
\begin{equation}
|\operatorname{arg}z| < \left(1 + \frac{2-a}{\nu}\right) \frac{\pi}{2} 
\end{equation}
as long as $1 + (2-a)/\nu > 0$, which is guaranteed for $a \leq 2$. 
\end{proof}

\medskip

As seen above, the contour $\mathcal{C}_{ic\infty}$ can be 
deformed into 
$\mathcal{C}_{+\infty}$ if $\nu < a$ or into $\mathcal{C}_{-\infty}$ if $\nu > a$. We will use this in conjunction with the residue theorem. We expect the same problem seen in the previous proof to occur here, that is, having poles and zeros of the integrand at the same points. We also expect to be able to use the eq.\eqref{ap.B.general.G.tau.2} to simplify the integrand. As a matter of fact, using eq.\eqref{ap.B.general.G.tau.2} we obtain 
\begin{equation}
\label{laplace.expression}
\begin{aligned}
& \mathcal{L}_t[ E_{a,m,l}\left(-\lambda t^\nu\right);z] =
\frac{G(\varphi+a\tau;\tau)}{G(\varphi;\tau)} \,\frac{1}{z} \,\frac{1}{2\pi i} 
\int_{\mathcal{C}} \left(\frac{z^\nu}{\lambda}\right)^{-s} 
\frac{\Gamma(1+\nu s)\Gamma(1+s) \Gamma(-s) G(\varphi+s;\tau)}{G(\varphi+a\tau+s;\tau)} ds. 
\end{aligned}
\end{equation}
The integrand has poles due to the gamma functions in the 
numerator located at $\acute{\mys}_n = n$, $\grave{\mys}^{(1)}_{n} = -(n+1)$, $\grave{\mys}^{(2)}_n = -(n+1)/\nu$ and poles due to the double gamma function in the denominator located at $\grave{\mys}^{(3)}_{m,n} = 
-(\varphi+a\tau+m\tau+n)$. The simplest case is clearly to use the contour $\mathcal{C}_{+\infty}$, and after evaluating the residues of the integrand at $\acute{\mys}_n = n$ we obtain the series in eq.\eqref{laplace.1}. As we saw, this can be done when $\nu < a$, in accordance with what we have  discussed about the convergence of that series. 

When $\nu > a$ we have to use the contour $\mathcal{C}_{-\infty}$ and then we have to evaluate the residues at $\grave{\mys}^{(1)}_{n} = -(n+1)$, $\grave{\mys}^{(2)}_n = -(n+1)/\nu$ and $\grave{\mys}^{(3)}_{m,n} = 
-(\varphi+a\tau+m\tau+n)$. The problem is that the order of the residues  depends on $\nu > 0$ and $\tau > 0$. 
In fact, if $\nu$ is rational, there will be poles $\grave{\mys}^{(1)}_n$ and $\grave{\mys}^{(2)}_{n^\prime}$ that will coincide, being therefore poles of order 2. If $\nu$ is irrational, the poles $\grave{\mys}^{(1)}_n$ and $\grave{\mys}^{(2)}_{n^\prime}$ will not coincide. But an even bigger problem is related to the poles $\grave{\mys}^{(3)}_{m,n}$. Besides the possibility of them coinciding with $ \grave{\mys}^{(1)}_{n^\prime}$ or $\grave{\mys}^{(2)}_{n^{\prime\prime}}$ depending on the value of $\varphi + a\tau$, another problem lies in the fact that the order of the poles $\grave{\mys}^{(3)}_{m,n}$ depends on $\tau$. If $\tau$ is irrational, these poles are simple, but if, for example, $\tau = 1$, the order of $\grave{\mys}^{(3)}_{m,n}$ is $(m+n+1)$.  

Notwithstanding, the first residue is not difficult to evaluate. As a concrete example, let us consider the case studied in \cite{BLV}, where $a=\alpha$, $m = 1+\gamma/\alpha$, $l = \gamma/\alpha$, $\nu = \alpha + \gamma$ with 
$0 < \alpha < 1$ and $\nu = \alpha + \gamma > 0$. We have 
\begin{equation}
\label{residue.first}
\begin{aligned}
& \mathcal{L}_t[ E_{\alpha,1+\frac{\gamma}{\alpha},\frac{\gamma}{\alpha}}\left(- \lambda t^{\nu}\right);z] \\[1ex]
& = 
\begin{cases}
\displaystyle \frac{\Gamma(1-\nu)}{\lambda \Gamma(1-\alpha)} \,\frac{1}{z^{1-\nu}} + 
\mathcal{O}(z^{\operatorname{min}(2\nu-1,0)}) , & \quad \; \alpha < \nu < 1 ,\\[2ex]
\displaystyle \frac{1}{\lambda (1-\alpha)}\left[-\log{z} + \log{\lambda} + 
\psi(1) + \alpha\left(1-\psi(1-\alpha)\right)\right] + \mathcal{O}(z) , & \quad \; 
\nu = 1 ,\\[2ex]
\displaystyle \frac{\pi/\nu}{\lambda^{1/\nu}\tau^{\alpha/\nu}\Gamma\left(\gamma/\nu\right) 
\sin\pi/\nu} + \mathcal{O}(z^{\operatorname{min}(1,\nu-1)}) , &\quad \; \nu > 1 ,
\end{cases}
\end{aligned}
\end{equation}
where $\psi(\cdot)$ is the digamma function. 

\begin{figure}[hbt]
\begin{center}
\includegraphics[width=12cm]{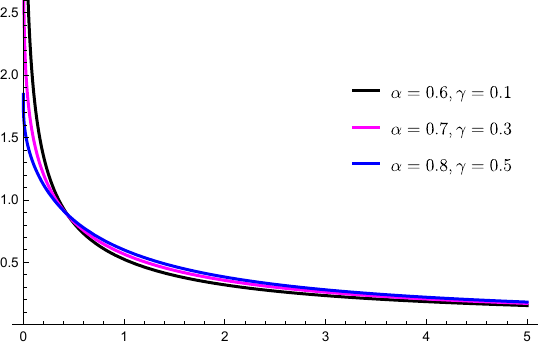}
\caption{Plots of the Laplace transform $\mathcal{L}_t[ E_{a,m,l}\left(-\lambda t^\nu\right);z] $ obtained from eq.\eqref{laplace.expression} for $\lambda=1$ 
and $\alpha$ and $\gamma$ as above. \label{fig.1}}
\end{center}
\end{figure}

 Despite the difficulty in obtaining a series representation for the Laplace transform of the KS function when $\nu > a$, if we are interested in using this representation to compute numerical values of the function, the best option is actually to use the contour integral itself in the complex plane. In fact, integral representations are frequently the most useful form for numerical computation \cite{Luke, Olver, Temme}, and we already have this expression. 
 In Figure~\ref{fig.1} are the plots obtained using 
 Mathematica 14.2 by numerically integrating eq.\eqref{laplace.expression} with $\lambda = 1$ and eq.\eqref{double.gamma.int.rep} for the values of $G(\cdot,\tau)$  
for the cases (i) $\alpha = 0.6$ and $\gamma = 0.1$ ($\alpha< \nu < 1$), (ii) $\alpha = 0.7$ and $\gamma = 0.3$ ($\nu = 1$) and (iii) $\alpha = 0.8$ and $\gamma = 0.5$ ($\nu > 1$). 

Note that the value of $\mathcal{L}_t[ E_{a,m,l}\left(-\lambda t^\nu\right);z]$ 
at $ z = 0$ is finite only when $\nu > 1$. From eq.\eqref{residue.first} we see that this value for $\alpha = 0.8$ and $\gamma = 0.5$ is $1.85457$.

\section{Conclusions}

We defined the Fox-Barnes $J$-function as a contour integral in the complex plane where the integrand is given by a ratio of products of gamma double gamma function involving several parameters, and showed how it generalizes the Fox $H$-function, which appears as an example of a completely hyper-balanced Fox-Barnes $J$-function. We studied the conditions for the existence of this contour integral. Some properties of the Fox-Barnes $J$-function were proved. As an application, we showed how the Laplace transform of the Kilbas-Saigo function can be conveniently written in terms of the Fox-Barnes $J$-function even in cases where the usual series representation is not convergent. As part of future studies, we hope to find new applications of the Fox-Barnes $J$-function, in particular in situations involving  variations in the definition of fractional derivatives \cite{BLV}.

\medskip

\bigskip
\bigskip

\noindent \textbf{\sffamily Acknowledgements:} We are grateful to FAPESP (Funda\c{c}\~ao de Amparo \`a Pesquisa do Estado de S\~ao Paulo) for the financial support (process 24/17510-6) and to Prof. L. Beghin and Sapienza Universit\`a di Roma for the hospitality during our stay as visiting professor.

\end{document}